\documentclass{cmam}
\usepackage{amsmath}
\usepackage{amssymb}
\usepackage{graphicx}

\issueinfo{12}{1}{2012}
\begin{document}

\title{EXPONENTIALLY CONVERGENT FUNCTIONAL-DISCRETE METHOD FOR EIGENVALUE TRANSMISSION PROBLEMS WITH DISCONTINUOUS FLUX AND POTENTIAL AS A FUNCTION IN SPACE  $L_1$}
\titlerunning{\CMAM{} Exponentially convergent functional-discrete method\ldots}
\author{V.L. Makarov \and N.O. Rossokhata \and D.V. Dragunov}
\institute{V.L. Makarov \at
\emph{Institute of Mathematics, National Academy of Sciences, 3 Tereschenkivska St., 01601 Kyiv, Ukraine}\\
\email{makarov@imath.kiev.ua}.
\and
N.O. Rossokhata \at
\emph{Department of Mathematics and Statistics,
Concordia University,
1455 De Maisonneuve Blvd. West, QC H3G 1M8, Canada}\\
\email{nataross@gmail.com}.
\and
D.V. Dragunov \at
\emph{Institute of Mathematics, National Academy of Sciences, 3 Tereschenkivska St., 01601 Kyiv, Ukraine} \\
\email{dragunovdenis@gmail.com}.
}
\begin{abstract}
  Based on the  functional-discrete technique (FD-method), an algorithm for eigenvalue transmission problems with discontinuous flux and integrable potential is developed. The case of the potential as a function belonging to the functional space $L_1$ is studied for both linear and nonlinear eigenvalue problems. The sufficient conditions providing  superexponential convergence rate of the method were obtained. Numerical examples are presented to support the theory. Based on the numerical examples and the convergence results, conclusion about analytical properties of eigensolutions for nonself-adjoint differential operators is made.
\end{abstract}
\maketitle

\section{Introduction}\label{s_1}

Physical and chemical processes in multilayer medium can be modeled by partial differential equations. Because of different physical and chemical properties of layer materials, besides initial and boundary conditions such models involve transmission (matching) conditions describing processes in the boundary of two layers. As a rule such conditions provide jump of a solution or flux, or both of them. Transmission problem for elliptic equations was proposed by M.M. Picone in 1954 \cite{picone} for the case of discontinuous solutions and continuous flux and then it was studied in \cite{campanato,lions, stampacchia, schechter}. Linear eigenvalue transmission problems were investigated in \cite{ozisik1, ozisik2, ozisik3}. The numerical treatment of transmission problems are given in \cite{mrb1, mrb2, ross4.2, mrrew, ross1, ross3.3, givoli, vulkov, jovan1, jovan2} et al. Particularly, papers \cite{mrb1, mrb2, ross4.2, mrrew, ross1, ross3.3} are devoted to the development of numerical functional-discrete exponentially convergent methods for the eigenvalue transmission problem with discontinuous solutions and continuous flux.

In this article we extend the FD-approach to the eigenvalue transmission problems with continuous solutions and discontinuous flux.

The paper is organized as follows. Section \ref{s_2} is devoted to the eigenvalue transmission problem with a potential belonging to the functional space $L_1(0,1)$. In section \ref{s_2} we formulate the problem with the additional differential condition providing a uniqueness of the solution, then we describe the numerical technique and prove the convergence theorem. In section \ref{S_3}  the corresponding nonlinear eigenvalue problem has been studied and  the convergence result similar to that of section \ref{s_2} has been obtained. The results of numerical experiments are discussed in section \ref{s_4}. They are in good agreement with the theoretical findings. Section \ref{s_5} includes conclusions about the proposed approach.


\section{Eigenvalue transmission problem with a linear potential as a function in $L_1$}\label{s_2}
\subsection{Formulation of the problem}

Let us consider the following eigenvalue problem
$$u_1''(x)+[\lambda-q(x)]u_1(x)=0, \hspace*{5mm} x\in \left(0,\frac{1}{2}\right), \hspace*{5mm} q(x)\in L_1(0,1)$$
\begin{equation}\label{3}
u_2''(x)+[\lambda-q(x)]u_2(x)=0, \hspace*{5mm} x\in \left(\frac{1}{2},1\right), \hspace*{22mm}
\end{equation}
with transmission conditions
$$\left[u\left(\frac{1}{2}\right)\right]=0, \left[u'\left(\frac{1}{2}\right)\right]=1,$$
and boundary conditions
$$u_1(0)=u_2(1)=0, u'_1(0)=1,$$
where $[f(x_0)]=f_2(x_0)-f_1(x_0)$ is a jump of a function $f$ at the point $x_0$.

\subsection{Description of the FD-method for the linear case}

According to the FD-approach \cite{makarov1} instead of the original problem (\ref{3}) we consider the following problem with parameter $t\in [0,1]$:
$$\frac{\partial^2}{\partial x^2}u_1(x,t)+[\lambda(t)-tq(x)]u_1(x,t)=0, \hspace*{5mm} x\in \left(0,\frac{1}{2}\right),$$
\begin{equation}\label{4}
\frac{\partial^2}{\partial x^2}u_2(x,t)+[\lambda(t)-tq(x)]u_2(x,t)=0, \hspace*{5mm} x\in \left(\frac{1}{2},1\right),
\end{equation}
$$u_1(0,t)=u_2(1,t)=0, \left.\frac{\partial}{\partial x}u^{\prime}_{1}(x,t)\right|_{x=0}=1,$$
$$\left[u\left(\frac{1}{2},t\right)\right]=0, \left[u^{\prime}\left(\frac{1}{2},t\right)\right]=1, \forall t\in \left[0, 1\right],$$
where $q(x)\in L_1(0,1)$.


Taking into account that for $t=0$ we can easy find the exact solution $(\lambda(0), u_i(x,0), i=1,2)$ to problem (\ref{4}) and that for $t=1$ the solution $(\lambda(1), u_i(x,1), i=1,2)$ to problem (\ref{4}) coincides with that to problem (\ref{3}), we can suppose that solution to problem (\ref{4}) can be found in the form of power series with respect to $t$
\begin{equation}\label{5}
\lambda_n(t)=\sum_{j=0}^\infty\lambda_n^{(j)}t^j, \hspace*{5mm} u_{ni}(x,t)=\sum_{j=0}^\infty u_{ni}^{(j)}(x)t^j, \hspace*{5mm} i=1,2.
\end{equation}

Setting $t=1$, we obtain
\begin{equation}\label{6}
\lambda_n=\lambda_n(1)=\sum_{j=0}^\infty\lambda_n^{(j)}, \hspace*{5mm} u_{ni}(x)=u_{ni}(x,1)=\sum_{j=0}^\infty u_{ni}^{(j)}(x), \hspace*{5mm} i=1,2
\end{equation}
provided that these series converge.

Thus, we can represent the approximate solution to problem (\ref{3}) as the corresponding truncated series
\begin{equation}\label{7}
\stackrel{m}{\lambda}_n=\sum_{j=0}^m\lambda_n^{(j)}, \hspace*{5mm} \stackrel{m}{u}_{ni}(x)=\sum_{j=0}^m u_{ni}^{(j)}(x), \hspace*{5mm} i=1,2,
\end{equation}
which are referred to as the approximations of the $m$-th rank.

Substituting \eqref{5} into \eqref{4} and equating the coefficients of equal powers of $t$, we obtain the following recursive sequence of problems:

for $j=-1$
$$u_{n1}^{(0)''}(x)+\lambda_n^{(0)}u_{n1}^{(0)}(x)=0, \hspace*{5mm} x\in \left(0,\frac{1}{2}\right),$$
\begin{equation}\label{8}
u_{n2}^{(0)''}(x)+\lambda_n^{(0)}u_{n2}^{(0)}(x)=0, \hspace*{5mm} x\in \left(\frac{1}{2},1\right),
\end{equation}
$$u_{n1}^{(0)}(0)=u_{n2}^{(0)}(1)=0, u_{n1}^{(0)'}(0)=1,$$
$$\left[u_n^{(0)}\left(\frac{1}{2}\right)\right]=0, \left[u_n^{(0)'}\left(\frac{1}{2}\right)\right]=1;$$
and
for $j=0,1,2,\ldots$
$$u_{n1}^{(j+1)''}(x)+\lambda_n^{(0)}u_{n1}^{(j+1)}(x)=F_{n1}^{(j+1)}(x,q), \hspace*{5mm} x\in \left(0,\frac{1}{2}\right),$$
\begin{equation}\label{9}
u_{n2}^{(j+1)''}(x)+\lambda_n^{(0)}u_{n2}^{(j+1)}(x)=F_{n2}^{(j+1)}(x,q), \hspace*{5mm} x\in \left(\frac{1}{2},1\right),
\end{equation}
$$u_{n1}^{(j+1)}(0)=u_{n2}^{(j+1)}(1)=0, u_{n1}^{(j+1)\prime}(0)=0,$$
$$\left[u_n^{(j+1)}\left(\frac{1}{2}\right)\right]=0, \left[u_n^{(j+1)\prime}\left(\frac{1}{2}\right)\right]=0,$$
where
\begin{equation}\label{9*}
F_{ni}^{(j+1)}(x,q)=-\sum_{p=0}^j\lambda_n^{(j+1-p)}u_{ni}^{(p)}(x)+q(x)u_{ni}^{(j)}(x), \hspace*{5mm} i=1,2.
\end{equation}

We obviously can find the exact solution $(\lambda_n^{(j)}, u_{ni}^{(j)}(x), i=1,2)$, $(j=0,1,\ldots)$ to each of the problems.

Thus, to apply the proposed algorithm one should find the rough approximation \newline$(\lambda_n^{(0)}, u_{ni}^{(0)}(x), i=1,2)$ as the exact solution to problem (\ref{8}) and then consecutively find the  corrections $(\lambda_n^{(j)}, u_{ni}^{(j)}(x), i=1,2)$ as the exact solutions to problems (\ref{9}), (\ref{9*}) for $j=0,1,2\ldots$.


\subsection{Convergence result}

Formulas (\ref{6}) and (\ref{7}) imply that the algorithm's error can be estimated in the following way
\begin{eqnarray}\label{series_of_norms}
    \|u_n-\stackrel{m}{u}_n\|\leq \sum_{j=m+1}^\infty\|u_n^{(j)}\|,\nonumber\\
    \|\lambda_n-\stackrel{m}{\lambda}_n\|\leq \sum_{j=m+1}^\infty|\lambda_n^{(j)}|,
\end{eqnarray}
provided that series  in the right sides are convergent. Therefore, our aim is to find conditions providing the convergence of the corresponding series and to estimate their convergence rates. To achieve that we will try to construct a geometric sequence with  denominator less than 1, which dominates the terms of series in formulas  \eqref{series_of_norms}.

Firstly, let us find the solution to problem \eqref{8}. Taking into account the differential equation, the boundary and matching, we can write eigenfunction of problem \eqref{8} in the following form
$$u_{n1}^{(0)}(x)=\frac{\sin\left(\sqrt{\lambda_n^{(0)}}x\right)}{\sqrt{\lambda_n^{(0)}}}, \hspace*{5mm} u_{n2}^{(0)}(x)=c_2^{(0)}\sin\left(\sqrt{\lambda_n^{(0)}}(1-x)\right),$$
where the unknown constant $c_2^{(0)}$ and eigenvalue $\lambda_n^{(0)}$ can be determined using the matching conditions
\begin{equation}\label{10}
-c_2^{(0)}\sin\left(\sqrt{\lambda_n^{(0)}}/2\right)+\frac{\sin\left(\sqrt{\lambda_n^{(0)}}/2\right)}{\sqrt{\lambda_n^{(0)}}}=0,
\end{equation}
$$c_2^{(0)}\sqrt{\lambda_n^{(0)}}\cos\left(\sqrt{\lambda_n^{(0)}}/2\right)+\cos\left(\sqrt{\lambda_n^{(0)}}/2\right)=-1.$$

From \eqref{10} it follows that we have two different sequences of eigensolutions of problem \eqref{8}:


\begin{enumerate}
\item[I.]
\begin{equation}\label{11}
\lambda_n^{(0)}=4\pi^2\left (\pm\frac{2}{3}+2n\right )^2, \hspace*{5mm} n=0,1,2\ldots
\end{equation}
$$u_{n1}^{(0)}(x)=\frac{\sin\left(\sqrt{\lambda_n^{(0)}}x\right)}{\sqrt{\lambda_n^{(0)}}}, \hspace*{5mm} u_{n2}^{(0)}(x)=\frac{\sin\left(\sqrt{\lambda_n^{(0)}}(1-x)\right)}{\sqrt{\lambda_n^{(0)}}}$$
and
\item[II.]
\begin{equation}\label{12}
\lambda_n^{(0)}=4\pi^2n^2,  \quad u_{n1}^{(0)}(x)=\frac{\sin\left(2\pi nx\right)}{2\pi n},
\end{equation}
$$  u_{n2}^{(0)}(x)=\frac{(-1)^{n+1}-1}{2\pi n}\sin\left(2\pi n(1-x)\right){=\frac{(-1)^{n}+1}{2\pi n}\sin\left(2\pi n x\right), \hspace*{5mm} n=1, 2\ldots}$$

\end{enumerate}

Next, we write the solution to the nonhomogeneous boundary value problem (\ref{9}) as
\begin{equation}\label{13}
u_{n1}^{(j+1)}(x)=\int\limits_0^x\frac{\sin\left(\sqrt{\lambda_n^{(0)}}(x-\xi)\right)}{\sqrt{\lambda_n^{(0)}}}F_{n1}^{(j+1)}(\xi)d\xi, \hspace*{5mm} x\in \left[0,\frac{1}{2}\right],
\end{equation}
$$u_{n2}^{(j+1)}(x)=c_2^{(j+1)}\sin\left(\sqrt{\lambda_n^{(0)}}(1-x)\right)- \hspace*{45mm}$$
$$- \int\limits_x^1\frac{\sin\left(\sqrt{\lambda_n^{(0)}}(x-\xi)\right)}{\sqrt{\lambda_n^{(0)}}}F_{n2}^{(j+1)}(\xi)d\xi, \hspace*{5mm} x\in \left[\frac{1}{2},1\right].$$

The unknown constant $c_2^{(j+1)}$ and correction $\lambda_n^{(j+1)}$ of the trial eigenvalue $\lambda_n$ can be found from the matching conditions of problem (\ref{9}).

Combining representations (\ref{13}) together with matching conditions (\ref{9}), we obtain the system of linear equations with respect to $c_{2}^{(j+1)}$ and $\lambda_n^{(j+1)}$
$$c_2^{(j+1)}\sin\left(\frac{\sqrt{\lambda_n^{(0)}}}{2}\right)=\hspace*{100mm}$$
$$=\int\limits_0^{1/2}\frac{\sin\left(\sqrt{\lambda_n^{(0)}}(\frac{1}{2}-\xi)\right)}{\sqrt{\lambda_n^{(0)}}}F_{n1}^{(j+1)}(\xi)d\xi +
\int\limits_{1/2}^1\frac{\sin\left(\sqrt{\lambda_n^{(0)}}(\frac{1}{2}-\xi)\right)}{\sqrt{\lambda_n^{(0)}}}F_{n2}^{(j+1)}(\xi)d\xi,$$
\begin{equation}\label{14}
c_2^{(j+1)}\cos\left(\frac{\sqrt{\lambda_n^{(0)}}}{2}\right)=\hspace*{80mm}
\end{equation}
$$=-\left[\int\limits_0^{1/2}\frac{\cos\left(\sqrt{\lambda_n^{(0)}}(\frac{1}{2}-\xi)\right)}{\sqrt{\lambda_n^{(0)}}}F_{n1}^{(j+1)}(\xi)d\xi +
\int\limits_{1/2}^1\frac{\cos\left(\sqrt{\lambda_n^{(0)}}(\frac{1}{2}-\xi)\right)}{\sqrt{\lambda_n^{(0)}}}F_{n2}^{(j+1)}(\xi)d\xi\right].$$

Now we are in position to find a recursive formula for the corrections $\lambda_n^{(j+1)}.$ For this purpose we have to consider several cases.

First, we consider the case when eigensolutions $(\lambda_n^{(0)}, u_{ni}^{(0)}(x), i=1,2)$ are determined by formula (\ref{11}). In this case we have that $\sin\left(\frac{\sqrt{\lambda_n^{(0)}}}{2}\right)\neq 0$ and $\cos\left(\frac{\sqrt{\lambda_n^{(0)}}}{2}\right)=-\dfrac{1}{2}\neq 0$. Hence, from system  \eqref{14} we get the following equation with respect to $\lambda^{(j+1)}_{n}:$
$$\int\limits_0^{1/2}\left [\frac{\sin\left(\sqrt{\lambda_n^{(0)}}(\frac{1}{2}-\xi)\right)}{\sin\left(\sqrt{\lambda_n^{(0)}}/2\right)}F_{n1}^{(j+1)}(\xi) + \frac{\cos\left(\sqrt{\lambda_n^{(0)}}(\frac{1}{2}-\xi)\right)}{\cos\left(\sqrt{\lambda_n^{(0)}}/2\right)}F_{n1}^{(j+1)}(\xi)\right ]d\xi+$$
$$+\int\limits_{1/2}^1\left [\frac{\sin\left(\sqrt{\lambda_n^{(0)}}(\frac{1}{2}-\xi)\right)}{\sin\left(\sqrt{\lambda_n^{(0)}}/2\right)}F_{n2}^{(j+1)}(\xi) + \frac{\cos\left(\sqrt{\lambda_n^{(0)}}(\frac{1}{2}-\xi)\right)}{\cos\left(\sqrt{\lambda_n^{(0)}}/2\right)}F_{n2}^{(j+1)}(\xi)\right ]d\xi =0,$$
which yields us the formula for $\lambda_{n}^{(j+1)}$
$$\lambda_n^{(j+1)}=\left\{\int\limits_0^{1/2}\frac{\sin\left(\sqrt{\lambda_n^{(0)}}x\right)\sin\left(\sqrt{\lambda_n^{(0)}} (1-x)\right)}{\sqrt{\lambda_n^{(0)}}}dx+ \int\limits_{1/2}^1\frac{\sin^2\left(\sqrt{\lambda_n^{(0)}}(1-x)\right)}{\sqrt{\lambda_n^{(0)}}}dx\right\}^{-1}\times$$
$$\times\left\{\int\limits_0^{1/2}\left [-\sum_{p=1}^j\lambda_n^{(j+1-p)}u_{n1}^{(p)}+q'(x)u_{n1}^{(j)}(x)\right ]\sin\left(\sqrt{\lambda_n^{(0)}}(1-x)\right)dx+\right .$$
$$\left .+\int\limits_{1/2}^1\left [-\sum_{p=1}^j\lambda_n^{(j+1-p)}u_{n2}^{(p)}+q'(x)u_{n2}^{(j)}(x)\right ]\sin\left(\sqrt{\lambda_n^{(0)}}(1-x)\right)dx\right\}.$$


Taking into account the equality
$$\int\limits_0^{1/2}\frac{\sin\left(\sqrt{\lambda_n^{(0)}}x\right)\sin\left(\sqrt{\lambda_n^{(0)}} (1-x)\right)}{\sqrt{\lambda_n^{(0)}}}dx+ \int\limits_{1/2}^1\frac{\sin^2\left(\sqrt{\lambda_n^{(0)}}(1-x)\right)}{\sqrt{\lambda_n^{(0)}}}dx=$$ $$=\frac{\left(\sin\left(\sqrt{\lambda^{(0)}_{n}}/2\right)\right)^{2}}{2\sqrt{\lambda^{(0)}_{n}}}= \frac{3}{8\sqrt{\lambda_n^{(0)}}},$$
we arrive to the following representation
\begin{equation}\label{15}
\lambda_n^{(j+1)}=\frac{8\sqrt{\lambda_n^{(0)}}}{3}
\left\{\int\limits_0^{1/2}\left [-\sum_{p=1}^j\lambda_n^{(j+1-p)}u_{n1}^{(p)}(x)+q(x)u_{n1}^{(j)}(x)\right ]\sin\left(\sqrt{\lambda_n^{(0)}}(1-x)\right)dx+\right .
\end{equation}
$$\left .+\int\limits_{1/2}^1\left [-\sum_{p=1}^j\lambda_n^{(j+1-p)}u_{n2}^{(p)}(x)+q(x)u_{n2}^{(j)}(x)\right ]\sin\left(\sqrt{\lambda_n^{(0)}}(1-x)\right)dx\right\}.$$

We can find unknown constant $c_2^{(j+1)}$ (see (\ref{13})) from either of equations (\ref{14}). For example, from the first one we have
\begin{equation}\label{15*}
c_2^{(j+1)}=\frac{1}{\sin(\sqrt{\lambda_n^{(0)}}/2)}\left (\int\limits_0^{1/2}\frac{\sin\left(\sqrt{\lambda_n^{(0)}}(\frac{1}{2}-\xi)\right)}{\sqrt{\lambda_n^{(0)}}} F_{n1}^{(j+1)}(\xi)d\xi +\right.$$
$$\left.+
\int\limits_{1/2}^1\frac{\sin\left(\sqrt{\lambda_n^{(0)}}(\frac{1}{2}-\xi)\right)}{\sqrt{\lambda_n^{(0)}}} F_{n2}^{(j+1)}(\xi)d\xi\right ).
\end{equation}
Formula \eqref{15*} yields us the estimate
\begin{eqnarray}\label{15**}
\left|c_2^{(j+1)}\right|\leq \frac{2}{\sqrt{3\lambda_n^{(0)}}}\left (\int\limits_0^{1/2}\left|F_{n1}^{(j+1)}(\xi)\right| d\xi + \int\limits_{1/2}^1\left|F_{n2}^{(j+1)}(\xi)\right| d\xi\right )=
\end{eqnarray}
$$ =\frac{2}{\sqrt{3\lambda_n^{(0)}}}\left (\left\|F^{(j+1)}_{n1}\right\|_{0,1(0,1/2)}+ \left\|F^{(j+1)}_{n2}\right\|_{0,1(1/2,1)}\right ) =
\frac{2}{\sqrt{3\lambda_n^{(0)}}}\left\|F^{(j+1)}_{n}\right\|_{0,1},$$
where $\|\cdot\|_{0,1(a,b)}$ denotes the norm in $L_1(a,b):$ $\|f(x)\|_{0,1,(a,b)}=\int_{a}^{b}|f(x)|d x;$ and $\|f\|_{0,1}=\|f_1\|_{0,1(0,1/2)}+\|f_2\|_{0,1(1/2,1)}$.

Hence, formulas (\ref{13}) - (\ref{15**}) imply the following  estimates
$$\|u_n^{(j+1)}\|_\infty =\max\left\{\|u_{n1}^{(j+1)}\|_{\infty, (0,1/2)}, \|u_{n2}^{(j+1)}\|_{\infty, (1/2,1)}\right\} \leq$$
$$\leq\max\left\{\frac{1}{\sqrt{\lambda_n^{(0)}}}\int_0^{1/2}|F_{n1}^{(j+1)}(x)|dx, |c_2^{(j+1)}|+\frac{1}{\sqrt{\lambda_n^{(0)}}}\int_{1/2}^1|F_{n2}^{(j+1)}(x)|dx\right\}\leq$$
$$\leq\frac{2}{\sqrt{3}\sqrt{\lambda_{n}^{(0)}}}\left\|F^{(j+1)}_{n}\left(x\right)\right\|_{0,1}+\frac{1}{\sqrt{\lambda_{n}^{(0)}}}\left\|F^{(j+1)}_{n}\left(x\right)\right\|_{0,1}\leq$$
$$\leq\frac{1}{\sqrt{\lambda_{n}^{(0)}}}\left(\frac{2}{\sqrt{3}}+1\right)\left\|F^{(j+1)}_{n}\left(x\right)\right\|_{0,1}\leq$$
$$\leq{\frac{2+\sqrt{3}}{\sqrt{3}}}\frac{1}{\sqrt{\lambda_n^{(0)}}}\left [\sum_{p=0}^j|\lambda_n^{(j+1-p)}|\cdot\|u_n^{(p)}\|_\infty +\|q\|_{0,1}\cdot\|u_n^{(j)}\|_\infty\right ],$$
$$|\lambda_n^{(j+1)}|\leq\frac{8}{3}\sqrt{\lambda_n^{(0)}}\left [\sum_{p=1}^j|\lambda_n^{(j+1-p)}\cdot\|u_n^{(p)}\|_\infty+\|q\|_{0,1}\cdot\|u_n^{(j)}\|_\infty\right ],$$
where $\|\cdot\|_{\infty,[a,b]}$ denotes the norm in $L_\infty[a,b]$ and
$\|f\|_{\infty}=\max\{\|f_1\|_{\infty,[0,1/2]}, \|f_2\|_{\infty,[1/2,1]}\}$.


Thus, to estimate $\|u_n^{(j+1)}\|_\infty$ and $|\lambda_n^{(j+1)}|$, we have to solve the following system of inequalities
\begin{equation}\label{16}
\|u_n^{(j+1)}\|_\infty\leq\frac{1}{a}\left [\sum_{p=0}^j|\lambda_n^{(j+1-p)}|\cdot\|u_n^{(p)}\|_\infty +\|q\|_{0,1}\cdot\|u_n^{(j)}\|_\infty\right ],
\end{equation}
$$|\lambda_n^{(j+1)}|\leq\frac{1}{b}\left [\sum_{p=1}^j\left|\lambda_n^{(j+1-p)}\right|\cdot\|u_n^{(p)}\|_\infty+\|q\|_{0,1}\cdot\|u_n^{(j)}\|_\infty\right ]$$
with $a=\frac{\sqrt{3}}{2+\sqrt{3}}\sqrt{\lambda_n^{(0)}}$, $b=\frac{3}{8}\frac{1}{\sqrt{\lambda_n^{(0)}}},$ $j=0,1,2\ldots.$

Introducing the new variables
\begin{equation}\label{17}
v_{j+1}=\frac{a^{j+1}}{b}\|u_n^{(j+1)}\|_\infty, \hspace*{5mm} \mbox{and} \hspace*{5mm} \mu_{j+1}=a^{j}|\lambda_n^{(j+1)}|
\end{equation}
we can rewrite system \eqref{16} in the following form
\begin{equation}\label{18}
v_{j+1}\leq\sum_{p=0}^j\mu_{j+1-p}v_p+\|q\|_{0,1}v_j,
\end{equation}
$$\mu_{j+1}\leq\sum_{p=1}^j\mu_{j+1-p}v_p+\|q\|_{0,1}v_j,\quad j=0,1,2\ldots.$$

Let us consider the scalar sequences $\{\overline{v}_{j}\}_{j=0}^{\infty},$ $\{\overline{\mu}_{j}\}_{j=1}^{\infty}$ defined by the following recurrence formulas
\begin{equation}\label{19}
\overline{v}_{j+1}=\sum_{p=0}^j\overline{\mu}_{j+1-p}\overline{v}_p+\|q\|_{0,1}\overline{v}_j,
\end{equation}
$$\overline{\mu}_{j+1}=\sum_{p=1}^j\overline{\mu}_{j+1-p}\overline{v}_p+\|q\|_{0,1}\overline{v}_j= \overline{v}_{j+1}-\overline{\mu}_{j+1}\overline{v}_0,\quad j=0,1,2,\ldots$$
with $\overline{v}_0=v_0=\|u^{(0)}_{n}/b\|_{\infty}=\frac{8}{3}$.

Comparing inequalities \eqref{18} with equalities \eqref{19}, it is easy to see that
$$v_{j+1}\leq \overline{v}_{j+1}, \hspace*{5mm}  \mu_{j+1}\leq \overline{\mu}_{j+1}.$$

Eliminating $\overline{\mu}_{j+1}$ from system (\ref{19}) we arrive at the recurrence formula
$$\overline{v}_{j+1}=\sum_{p=1}^j\overline{v}_{j+1-p}\overline{v}_p+\|q\|_{0,1}\overline{v}_j(1+\overline{v}_0),\quad j=0,1,2\ldots.$$

Multiplying both sides of the last equality by $z^{j+1}$  and summing them over $j$ from $0$ to $\infty$ we obtain the following equation
$$f(z)-\overline{v}_0=[f(z)-\overline{v}_0]^2+(1+\overline{v}_0)\|q\|_{0,1}zf(z)$$
or, in more convenient form,
$$[f(z)-\overline{v}_0]^2-[1-(1+\overline{v}_0)\|q'\|_{0,1}z]\cdot [f(z)-\overline{v}_0]+\overline{v}_0(1+\overline{v}_0)\|q'\|_{0,1}z=0,$$
where
\begin{equation}\label{20}
f(z)=\sum_{j=0}^\infty z^j\overline{v}_j.
\end{equation}
We have the quadratic equation with respect to $f(z)-\overline{v}_0$ with the roots
$$\left(f\left(z\right)-\overline{v}_{0}\right)_{1,2}=\frac{\left(1-\left(1+\overline{v}_{0}\right)\left\|q\right\|_{0,1}z\right)\pm \sqrt{D}}{2},$$
where
\begin{eqnarray}\label{D_1}
    D&=&\left(1+\overline{v}_{0}\right)^{2}\left\|q\right\|_{0,1}^{2}\times\\
    &&\times\left(\frac{1+2\overline{v}_{0}+2\sqrt{\overline{v}_{0}\left(1+\overline{v}_{0}\right)}} {\left(1+\overline{v}_{0}\right)\left\|q\right\|_{0,1}}-z\right)\left(\frac{1+2\overline{v}_{0} -2\sqrt{\overline{v}_{0}\left(1+\overline{v}_{0}\right)}}{\left(1+\overline{v}_{0}\right) \left\|q\right\|_{0,1}}-z\right).\nonumber
\end{eqnarray}
It is easy to see that the solution which represents the generating function \eqref{20} is the following one
\begin{equation}\label{D_2}
    f\left(z\right)=\overline{v}_{0}+\frac{\left(1-\left(1+\overline{v}_{0}\right)\left\|q\right\|_{0,1}z\right)- \sqrt{D}}{2}.
\end{equation}
It is obvious that the right-hand side of equality \eqref{D_2} can be expanded as a power series in $z$ where $z\in \left[0, R\right]$ and
\begin{eqnarray}\label{21}
    {R=\frac{1+2\overline{v}_{0}-2\sqrt{\overline{v}_{0}\left(1+\overline{v}_{0}\right)}}{\left(1+\overline{v}_{0}\right) \left\|q\right\|_{0,1}}=\frac{1}{\left(1+\overline{v}_{0}\right)\left\|q\right\|_{0,1}\left(1+2\overline{v}_{0}+ 2\sqrt{\overline{v}_{0}\left(1+\overline{v}_{0}\right)}\right)}=}\nonumber\\
    {=\frac{1}{\left(\frac{8}{3}+1\right)\left\|q\right\|_{0,1}\left(1+\frac{16}{3}+2\sqrt{\frac{88}{9}}\right)}.}
\end{eqnarray}
Here $R$ denotes the radius of convergence for series \eqref{20}.
On the other hand, it follows from convergence of series (\ref{20}) that $$R^{j}\overline{v}^j\leq\frac{c}{j^{1+\varepsilon}}$$
for some positive constants $c$ and $\varepsilon$.

Returning to (\ref{17}), we arrive at the estimates
\begin{equation}\label{u_estimate_1}
    \|u_n^{(j+1)}\|_\infty =\frac{b v_{j+1}}{a^{j+1}}\leq\frac{b}{a^{j+1}}\overline{v}_{j+1}\leq \frac{{3} c}{{8} (j+1)^{1+\varepsilon}\sqrt{\lambda_n^{(0)}}}\cdot\frac{1}{\left (\frac{\sqrt{3}\sqrt{\lambda_{n}^{(0)}}}{2+\sqrt{3}}R\right )^{j+1}}
\end{equation}
and
\begin{equation}\label{lambda_estimate_1}
   |\lambda_n^{(j+1)}|=\frac{\mu_{j+1}}{a^j}\leq\frac{\overline{\mu}_{j+1}}{a^j}\leq \frac{\overline{v}_{j+1}}{(1+\overline{v}_0)a^j}\leq \frac{c}{(j+1)^{1+\varepsilon}R}\cdot\frac{1}{\left (\frac{\sqrt{3}\sqrt{\lambda_{n}^{(0)}}}{2+\sqrt{3}}R\right )^j}.
\end{equation}
Thus, series (\ref{6}) converge, and the last two inequalities imply the error estimates
$$\|u_n-\stackrel{m}{u}_n\|_\infty \leq\sum_{j=m+1}^\infty\|u_n^{(j)}\|_{\infty}\leq \frac{{3} C}{{8} (m+1)^{1+\varepsilon}\sqrt{\lambda_n^{(0)}}}\cdot
\left (\frac{2+\sqrt{3}}{\sqrt{3}\sqrt{\lambda_{n}^{(0)}}R}\right )^{m+1},$$
\begin{equation}\label{22}
|\lambda_n-\stackrel{m}{\lambda_n}|\leq\sum_{j=m+1}^\infty |\lambda_n^{(j)}|\leq \frac{C}{(m+1)^{1+\varepsilon}R}\left(\frac{2+\sqrt{3}}{\sqrt{3}\sqrt{\lambda_{n}^{(0)}}R}\right)^m
\end{equation}
provided that
\begin{equation}\label{23}
r_n=\frac{2+\sqrt{3}}{\sqrt{3}\sqrt{\lambda_{n}^{(0)}}R}<1,
\end{equation}
where $$C=\frac{c}{1-\frac{2+\sqrt{3}}{\sqrt{3}\sqrt{\lambda_{n}^{(0)}}R}}.$$

Substituting expression for $R$ \eqref{21} into formula \eqref{23} we arrive at the following convergence condition for  the algorithm \eqref{7}, \eqref{8}, \eqref{9}
\begin{equation}\label{24}
{r_n=\frac{(2+\sqrt{3})\frac{11}{3}\left\|q \right\|_{0,1}\left(1+\frac{16}{3}+2\sqrt{\frac{88}{9}}\right)}{\sqrt{3}\sqrt{\lambda_{n}^{(0)}}}< 1.}
\end{equation}


Now let us consider the case when eigensolutions $(\lambda_n^{(0)}, u_{ni}^{(0)}(x), i=1,2)$ are determined by formulas (\ref{12}). We see that in this case $\sin\left(\sqrt{\lambda_n^{(0)}}/2\right)=0$. Hence, to find $\lambda_n^{(j+1)}$ we should use the first equation of (\ref{14}) and parameter $c_2^{(j+1)}$ can be determined from the second equation of (\ref{14}), that is,
\begin{equation}\label{25}
{\lambda_n^{(j+1)}=\left\{\int\limits_0^{1/2}\frac{\sin^{2}\left(2\pi nx\right)}{2\pi n}dx+((-1)^{n}+1) \int\limits_{1/2}^1\frac{\sin^{2}\left(2\pi n x\right)}{2\pi n}dx\right\}^{-1}\times}
\end{equation}
$${\times\left\{\int\limits_0^{1/2}\left [-\sum_{p=1}^j\lambda_n^{(j+1-p)}u_{n1}^{(p)}(x)+q(x)u_{n1}^{(j)}(x)\right ]\sin \left(2\pi n x\right)dx+\right .}$$
$${\left .+\int\limits_{1/2}^1\left [-\sum_{p=1}^j\lambda_n^{(j+1-p)}u_{n2}^{(p)}(x)+q(x)u_{n2}^{(j)}(x)\right ]\sin \left(2\pi n x\right)dx\right\},}$$
\begin{equation}\label{26}
{c_2^{(j+1)}=-\int_0^{1/2}\frac{\cos\left(2\pi n x\right)}{2\pi n}F_{n1}^{(j+1)}(x)dx -
\int_{1/2}^1\frac{\cos\left(2\pi n x\right)}{2\pi n}F_{n2}^{(j+1)}(x)dx.}
\end{equation}

It is easy to verify that
$$\int\limits_0^{1/2}\frac{\sin^{2}\left(2\pi nx\right)}{2\pi n}dx+((-1)^{n}+1) \int\limits_{1/2}^1\frac{\sin^{2}\left(2\pi n x\right)}{2\pi n}dx=\left\{
                                                                     \begin{array}{l}
                                                                       \cfrac{3}{8 \pi n},\quad n=2k,\;k\in \mathbb{N}, \\ \\
                                                                       \cfrac{1}{8 \pi n},\quad n=2k+1,\;k\in \mathbb{N}. \\
                                                                     \end{array}
                                                                   \right.
$$

If $n$ is an even number, then $u_{n2}^{(0)}(x)=\cfrac{\sin\left( 2\pi n x\right)}{\pi n}$. Substituting the last expression into formula (\ref{26}), we get
\begin{equation}\label{27}
\lambda_n^{(j+1)}=\frac{8\pi n}{3}
\left\{\int\limits_0^{1/2}\left [-\sum_{p=1}^j\lambda_n^{(j+1-p)}u_{n1}^{(p)}+q(x)u_{n1}^{(j)}(x)\right ]\sin \left(2\pi n x \right)dx+\right .
\end{equation}
$$\left .+\int\limits_{1/2}^1\left [-\sum_{p=1}^j\lambda_n^{(j+1-p)}u_{n2}^{(p)}+q(x)u_{n2}^{(j)}(x)\right ]\sin \left(2\pi n x\right)dx\right\}.$$
Formulas (\ref{13}), (\ref{26}) and (\ref{27}) imply that
$$\|u_n^{(j+1)}\|_\infty =\max\left\{\|u_{n1}^{(j+1)}\|_{\infty, (0,1/2)}, \|u_{n2}^{(j+1)}\|_{\infty, (1/2,1)}\right\} \leq$$
$$\leq\max\left\{\frac{1}{\sqrt{\lambda_n^{(0)}}}\int_0^{1/2}|F_{n1}^{(j+1)}(x)|dx, |c_2^{(j+1)}|+\frac{1}{\sqrt{\lambda_n^{(0)}}}\int_{1/2}^1|F_{n2}^{(j+1)}(x)|dx\right\}\leq$$
\begin{equation}\label{D_10_11_10}
    \leq \left(\frac{1}{2\pi n}+\frac{1}{2 \pi n}\right)\left\|F^{(j+1)}_{n}\left(x\right)\right\|_{0,1}\leq
\end{equation}
$$\leq {\frac{1}{\pi n}}\left [\sum_{p=0}^j|\lambda_n^{(j+1-p)}|\cdot \|u_n^{(p)}\|_\infty + \|q\|_{0,1}\|u_n^{(j)}\|_\infty\right ]$$
and
\begin{equation}\label{D_10_11_11}
    |\lambda_n^{(j+1)}|\leq \frac{8\pi n}{3}\left [\sum_{p=1}^j|\lambda_n^{(j+1-p)}|\cdot \|u_n^{(p)}\|_\infty + \|q\|_{0,1}\|u_n^{(j)}\|_\infty\right ],
\end{equation}
that is, we get inequalities (\ref{16}) with
$$a=\pi n, \hspace*{5mm} b=\frac{3}{8\pi n}.$$


Repeating the same computations as for the previous case, we obtain the similar error estimates
$$\|u_n-\stackrel{m}{u}_n\|_\infty\leq \frac{{3} C}{{8 \pi n} (m+1)^{1+\varepsilon}}
\left (\frac{1}{ \pi nR}\right )^{m+1},$$
\begin{equation}\label{28}
|\lambda_n-\stackrel{m}{\lambda_n}|\leq \frac{C}{(m+1)^{1+\varepsilon}R}\left (\frac{1}{\pi nR}\right )^m
\end{equation}
provided that
\begin{equation}\label{29}
r_n=\frac{1}{\pi nR}=\frac{\frac{11}{3}\left\|q\right\|_{0,1}\left(1+\frac{16}{3}+2\sqrt{\frac{88}{9}}\right)}{\pi n}< 1,
\end{equation}
where
\begin{equation}\label{D_29'}
    C=\frac{c}{1-\frac{1}{\pi nR}}.
\end{equation}

If $n$ is an odd number, then inequality \eqref{D_10_11_11} should be replaced by the following one
\begin{equation}\label{D_10_11_12}
    |\lambda_n^{(j+1)}|\leq 8\pi n\left [\sum_{p=1}^j|\lambda_n^{(j+1-p)}|\cdot \|u_n^{(p)}\|_\infty + \|q\|_{0,1}\|u_n^{(j)}\|_\infty\right ],
\end{equation}
however, the estimate for $\|u_n^{(j+1)}\|_\infty$ remains the same, see \eqref{D_10_11_10}.
Thus we arrive at the equation (\ref{16}) with
$$a=\pi n, \hspace*{5mm} b=\frac{1}{8\pi n}$$
and we get the following error estimates
$$\|u_n-\stackrel{m}{u}_n\|_\infty\leq \frac{C}{{ 8\pi n} (m+1)^{1+\varepsilon}}
\left ({\frac{1}{ \pi nR}}\right )^{m+1},$$
\begin{equation}\label{31}
|\lambda_n-\stackrel{m}{\lambda_n}|\leq \frac{C}{(m+1)^{1+\varepsilon}R}\left ({\frac{1}{ \pi nR}}\right )^m
\end{equation}
with $C$ defined in \eqref{D_29'},
provided that inequality \eqref{29} holds true.

Therefore, we have obtained the following convergence result.\\
\begin{theorem}
\indent
\begin{enumerate}

\item[(a)] Assume that the index $n$ of a trial eigenpair satisfies condition (\ref{24}) with $\lambda_n^{(0)}=4\pi^2(\pm\frac{2}{3}+2n)^2$, $(n=0,1,\ldots).$ Then the numerical algorithm (\ref{7}), (\ref{11}), (\ref{13}), (\ref{9*}), (\ref{15}), (\ref{15*}) converges to the corresponding eigensolution of problem (\ref{3}) super-exponentially with error estimate (\ref{22}).

\item[(b)] Assume that the index $n$ of a trial eigenpair satisfies condition (\ref{29}) with $\lambda_n^{(0)}=4\pi^2n^2,$ $n=1,2,\ldots.$ Then the numerical algorithm (\ref{7}), (\ref{12}), (\ref{13}), (\ref{9*}), (\ref{25}), (\ref{26}) converges to the corresponding eigensolution of problem (\ref{3}) super-exponentially with error estimate (\ref{28}) for even number $n$ and with error estimate (\ref{31}) for odd number $n$.
\end{enumerate}
\end{theorem}

     Using the term ``super-exponential convergence'', we intend to emphasize the fact that owing to the presence of factor $(m+1)^{1+\varepsilon}$ in the denominators of estimates \eqref{22}, \eqref{28} and \eqref{31} the method remains convergent even if $r_{n}=1.$

\section{Eigenvalue transmission problem with a nonlinear potential including a function in space $L_1$}\label{S_3}
\subsection{Formulation of the problem}

Let us consider the following eigenvalue problem

$$u_1''(x)+[\lambda-q(x)]u_1(x)-N(u_1(x))=0, \hspace*{5mm} x\in \left(0,\frac{1}{2}\right), \hspace*{5mm} q(x)\in L_1(0,1)$$
\begin{equation}\label{33}
u_2''(x)+[\lambda-q(x)]u_2(x)-N(u_2(x))=0, \hspace*{5mm} x\in \left(\frac{1}{2},1\right), \hspace*{19mm}
\end{equation}
$$u_1(0)=u_2(1)=0, u'_1(0)=1, \hspace*{45mm}$$
$$\left[u\left(\frac{1}{2}\right)\right]=0, \left[u'\left(\frac{1}{2}\right)\right]=1, \hspace*{45mm}$$
where $N: \mathbb{R} \to \mathbb{R} $ is an analytical function with respect to $u$, $N(0)=0,$ that is,
\begin{equation}\label{34}
{N(u)=\sum\limits_{i=1}^{\infty}a_{i}u^{i},\quad \forall u\in \mathbb{R},\quad a_{i}\in \mathbb{R},\quad  \forall i\in \mathbb{N}.}
\end{equation}
By $\overline{N}\left(u\right)$ we denote a fixed function from $C^{\infty}\left(\mathbb{R}\right)$ such that
\begin{equation}\label{35'}
   \left\|\frac{d^{n}}{d u^{n}}N\left(u\right)\right\|\leq \left.\frac{d^{n}}{d u^{n}}\overline{N}\left(u\right)\right|_{u=\left|u\right|},\quad n=0,1,\ldots.
\end{equation}
Such a function exists. Moreover, as it follows from \eqref{34}, the function $\overline{N}\left(u\right)=\sum\limits_{i=1}^{\infty}\left|a_{i}\right|u^{i}$ satisfies inequalities \eqref{35'} for any $n\in \mathbb{N}$.


\subsection{Description of the FD-method for the nonlinear case}

Using the general idea of FD-method, we approximate a numerical solution to problem \eqref{33} by the truncated series \eqref{7} where an eigenpair $(\lambda_n^{(0)},u_{ni}^{(0)}(x), i=1,2)$ is the solution to problem \eqref{8} (see formulas \eqref{11} and  \eqref{12}) and pairs $(\lambda_n^{(j+1)},u_{ni}^{(j+1)}(x), i=1,2)$ $(j=0,1,2,\ldots )$ are the solutions to the nonhomogeneous transmission problems (\ref{9}) with
\begin{equation}\label{37}
F_{ni}^{(j+1)}(x,q)=-\sum_{p=0}^j\lambda_n^{(j+1-p)}u_{ni}^{(p)}(x)+q(x)u_{ni}^{(j)}(x)+ A_{j}(N; u_{ni}^{(0)},...,u_{ni}^{(j)}),
\end{equation}
where
\begin{equation}\label{38}
A_{j}(N; u_{ni}^{(0)},...,u_{ni}^{(j)})=
\end{equation}
$$=\sum\limits_{\substack{\alpha_1+\ldots +\alpha_j=j \\ \alpha_1\geq \ldots \geq \alpha_j\geq 0 \\ \alpha_i\in \mathbb{Z}}} N_u^{(\alpha_1)}(u_{ni}^{(0)}(x))\frac{[u_{ni}^{(1)}(x)]^{\alpha_1-\alpha_2}}{(\alpha_1-\alpha_2)!} \cdots
\frac{[u_{ni}^{(j-1)}(x)]^{\alpha_{j-1}-\alpha_{j}}}{(\alpha_{j-1}-\alpha_{j})!}
\cdot \frac{[u_{ni}^{(j)}(x)]^{\alpha_j}}{(\alpha_j)!}, \hspace*{3mm} j>0,$$
$$A_{0}(u_{ni}^{(0)})=N(u_{ni}^{(0)})$$
are the Adomian polynomials (see \cite{seng1}, \cite{seng2}).

\subsection{Convergence result}

As it was in the linear case described in section \ref{s_2}, zero approximation $(\lambda_n^{(0)},u_{ni}^{(0)}(x), i=1,2)$ can be found using formulas (\ref{11}) and (\ref{12}). Applying the same approach as in section \ref{s_2}, we determine $u_{ni}^{(j+1)}(x)$ $(i=1,2)$ according to formulas (\ref{13}), (\ref{37}), (\ref{38}). The unknown constants $c_2^{(j+1)}$ and $\lambda_n^{(j+1)}$ can be found using the matching conditions (\ref{14}), (\ref{37}), (\ref{38}).

Hence, for the case of eigenpairs $(\lambda_n^{(0)},u_{ni}^{(0)}(x), i=1,2)$ determined by formulas \eqref{11}, we obtain the following expression for $\lambda_n^{(j+1)}$
\begin{equation}\label{40}
\lambda_n^{(j+1)}=\frac{8\sqrt{\lambda_n^{(0)}}}{3}\times \hspace*{80mm}
\end{equation}
$$\times\left\{\int_0^{1/2}\left [-\sum_{p=1}^j\lambda_n^{(j+1-p)}u_{n1}^{(p)}+q(x)u_{n1}^{(j)}(x)+ A_{j}(N; u_{n1}^{(0)},...,u_{n1}^{(j)})\right ]\sin\sqrt{\lambda_n^{(0)}}(1-x)dx+\right .$$
$$+\left .\int_{1/2}^1\left [-\sum_{p=1}^j\lambda_n^{(j+1-p)}u_{n2}^{(p)}+q(x)u_{n2}^{(j)}(x)+ A_{j}(N; u_{n2}^{(0)},...,u_{n2}^{(j)})\right ]\sin\sqrt{\lambda_n^{(0)}}(1-x)dx\right\},$$
and parameter $c_2^{(j+1)}$ is determined by formula \eqref{15*} with \eqref{37}.

Similarly to the linear case discussed in section \ref{s_2},  formulas (\ref{13}), (\ref{15*}), (\ref{37}) together with (\ref{34}) and (\ref{38}) yield us the following system of inequalities
$$
\|u_n^{(j+1)}\|_\infty\leq\frac{1}{a}\Biggl [\sum_{p=0}^j|\lambda_n^{(j+1-p)}|\cdot\|u_n^{(p)}\|_\infty +\|q\|_{0,1}\cdot\|u_n^{(j)}\|_\infty+\Biggr.$$
\begin{equation}\label{41}
\Biggl . +\sum\limits_{\substack{\alpha_1+\ldots +\alpha_j=j \\ \alpha_1\geq \ldots \geq \alpha_j\geq 0 \\ \alpha_i\in \mathbb{Z}}}N_u^{(\alpha_1)}(\|u_n^{(0)}\|_\infty ) \frac{\|u_n^{(1)}\|_\infty^{\alpha_1-\alpha_2}}{(\alpha_1-\alpha_2)!} \cdots
\frac{\|u_n^{(j-1)}\|_\infty^{\alpha_{j-1}-\alpha_{j}}}{(\alpha_{j-1}-\alpha_{j})!}
\cdot \frac{\|u_n^{(j)}\|_\infty^{\alpha_j}}{(\alpha_j)!}\Biggr],
\end{equation}
$$|\lambda_n^{(j+1)}|\leq\frac{1}{b}\Biggl [\sum_{p=1}^j|\lambda_n^{(j+1-p)}|\cdot\|u_n^{(p)}\|_\infty+\|q\|_{0,1}\cdot\|u_n^{(j)}\|_\infty+\Biggr.$$
$$\Biggl. +\sum\limits_{\substack{\alpha_1+\ldots +\alpha_j=j \\ \alpha_1\geq \ldots \geq \alpha_j\geq 0 \\ \alpha_i\in \mathbb{Z}}}N_u^{(\alpha_1)}(\|u_n^{(0)}\|_\infty ) \frac{\|u_n^{(1)}\|_\infty^{\alpha_1-\alpha_2}}{(\alpha_1-\alpha_2)!} \cdots
\frac{\|u_n^{(j-1)}\|_\infty^{\alpha_{j-1}-\alpha_{j}}}{(\alpha_{j-1}-\alpha_{j})!}
\cdot \frac{\|u_n^{(j)}\|_\infty^{\alpha_j}}{(\alpha_j)!}\Biggr]$$
with $a=\frac{\sqrt{3}}{2+\sqrt{3}}\sqrt{\lambda_n^{(0)}}$, $b=\frac{3}{8}\frac{1}{\sqrt{\lambda_n^{(0)}}}$.\\
Introducing new variables
\begin{equation}\label{42}
v_{j+1}=\frac{a^{j+1}}{b}\|u_n^{(j+1)}\|_\infty \quad  \mbox{and} \quad \mu_{j+1}=a^j|\lambda_n^{(j+1)}|,\quad j\geq-1,
\end{equation}
we arrive at the following system of inequalities:
\begin{equation}\label{43}
v_{j+1}\leq\sum_{p=0}^j\mu_{j+1-p}v_p+\|q\|_{0,1}v_j+A_{j}(\overline{N};v_0,v_1,...,v_j),
\end{equation}
$$\mu_{j+1}\leq\sum_{p=1}^j\mu_{j+1-p}v_p+\|q\|_{0,1}v_j+A_{j}(\overline{N};v_0,v_1,...,v_j),\quad j=1,2,\ldots.$$
To obtain similar estimates for $v_{1}$ and $\mu_{1}$ let us consider inequalities \eqref{41} with $j=0$ in more detail:
\begin{equation}\label{D_43}
  {  \|u_n^{(1)}\|_\infty \leq \frac{1}{a}\left[\left(|\lambda_n^{(1)}| +\left\|q\right\|_{0,1}\right)\left\|u_n^{(0)}\right\|_\infty+\overline{N}\left(\left\|u^{(0)}_{n}\right\|_{\infty}\right)\right],}
\end{equation}
\begin{equation}\label{D_44}
   { |\lambda_n^{(1)}|_\infty \leq \frac{1}{b}\left[\left\|q\right\|_{0,1}\left\|u_n^{(0)}\right\|_\infty+\overline{N}\left(\left\|u^{(0)}_{n}\right\|_{\infty}\right)\right].}
\end{equation}
Using notations \eqref{42} from \eqref{D_43} we obtain
\begin{eqnarray}\label{D_45}
  v_{1} &\leq & {\left(\mu_{1} +\left\|q\right\|_{0,1}\right)v_{0}+\frac{1}{b}\overline{N}\left(\|u^{(0)}_{n}\|_{\infty}\right)=}\nonumber
  \\ &=& \left(\mu_{1} +\left\|q\right\|_{0,1}\right)v_{0}+\frac{1}{b}\left(\overline{N}\left(\|u^{(0)}_{n}\|_{\infty}\right)-\overline{N}\left(0\right)\right)= \\
   &=& \left(\mu_{1} +\left\|q\right\|_{0,1}\right)v_{0}+\frac{1}{b}\overline{N}^{\prime}\left(\theta\|u^{(0)}_{n}\|_{\infty}\right)\|u^{(0)}_{n}\|_{\infty}\leq \nonumber\\
   &\leq & {\left(\mu_{1} +\left\|q\right\|_{0,1}\right)v_{0}+\overline{N}^{\prime}\left(v_{0}\right)v_{0}.} \nonumber
\end{eqnarray}
Similarly to \eqref{D_45} we can get the following inequality for $\mu_{1}$:
\begin{equation}\label{D_46}
  { \mu_{1}\leq \left\|q\right\|_{0,1}v_{0}+\overline{N}^{\prime}\left(v_{0}\right)v_{0}.}
\end{equation}

Now let us consider the sequences $\{\overline{v}_{j}\}_{j=0}^{\infty}$ and $\{\overline{\mu}_{j}\}_{j=0}^{\infty}$ defined by the following recurrence equalities
\begin{equation}\label{44}
\overline{v}_{j+1}=\sum_{p=0}^j\overline{\mu}_{j+1-p}\overline{v}_p+\|q\|_{0,1}\overline{v}_j+ A_j(\overline{N};\overline{v}_0,\overline{v}_1,...,\overline{v}_j),
\end{equation}
\begin{equation}\label{D_49}
    \overline{\mu}_{j+1}=\sum_{p=1}^j\overline{\mu}_{j+1-p}\overline{v}_p+\|q\|_{0,1}\overline{v}_j +A_j(\overline{N};\overline{v}_0,\overline{v}_1,...,\overline{v}_j)= \overline{v}_{j+1}-\overline{\mu}_{j+1}\overline{v}_0
\end{equation}
for $j\in \mathbb{N}$ and
\begin{equation}\label{D_47}
    \overline{v}_{1}=\left(\mu_{1} +\left\|q\right\|_{0,1}\right)\overline{v}_{0}+\overline{N}^{\prime}\left(\overline{v}_{0}\right)\overline{v}_{0},
\end{equation}
\begin{equation}\label{D_48}
    \overline{\mu}_{1}=\left\|q\right\|_{0,1}\overline{v}_{0}+\overline{N}^{\prime}\left(\overline{v}_{0}\right)\overline{v}_{0}=\overline{v}_{1}-\overline{\mu}_{1}\overline{v}_{0},
\end{equation}
where $\overline{v}_0=v_0=\frac{8}{3}$.
It is easy to see that $v_{j+1}\leq\overline{v}_{j+1}$ and $\mu_{j+1}\leq\overline{\mu}_{j+1}$ $\forall j\in \mathbb{N}\cup \left\{0\right\}.$

Eliminating variables $\overline{\mu}_{j+1}$ from equalities (\ref{44})-(\ref{D_48}) we arrive at the recurrence formulas
\begin{equation}\label{D_50}
  \overline{v}_{j+1}=\sum_{p=1}^j\overline{v}_{j+1-p}\overline{v}_p+(1+\overline{v}_0)\left(\|q\|_{0,1}\overline{v}_j+A_j(\overline{N};\overline{v}_0,\overline{v}_1,...,\overline{v}_j)\right)
\end{equation}
for $j\in \mathbb{N}$ and
 \begin{equation}\label{D_51}
    \overline{v}_{1}= \left(1+\overline{v}_{0}\right)\left(\left\|q\right\|_{0,1}\overline{v}_{0}+\overline{N}^{\prime} \left(\overline{v}_{0}\right)\overline{v}_{0}\right).
 \end{equation}

Applying the technique similar to that used in section 2 to equalities \eqref{D_50}, \eqref{D_51}
we can obtain the following nonlinear equation with respect to generating function $f(z)$
\begin{equation}\label{D_52}
    f(z)-\overline{v}_0=[f(z)-\overline{v}_0]^2+(1+\overline{v}_0)z[\|q\|_{0,1}f(z)+\overline{N}(f(z))+\overline{N}^{\prime}\left(\overline{v}_{0}\right)\overline{v}_{0}-\overline{N}\left(\overline{v}_{0}\right)],
\end{equation}
\begin{equation}\label{45}
f(z)=\sum_{j=0}^\infty z^j\overline{v}_j.
\end{equation}

Let us prove that power series \eqref{45} possesses the nonzero radius of convergence.  For this purpose we consider the inverse mapping $z=f^{-1}$. From equation \eqref{D_52} it follows that
\begin{equation}\label{D_53}
    z\left(f\right)= \frac{f-\overline{v}_0-[f-\overline{v}_0]^2}{(1+\overline{v}_0)[\|q\|_{0,1}f+\overline{N}(f)+\overline{N}^{\prime}\left(\overline{v}_{0}\right)\overline{v}_{0}-\overline{N}\left(\overline{v}_{0}\right)]}.
\end{equation}
Since function $z(f)$ \eqref{D_53} is holomorphic in some open interval containing the point $f_{0}=\overline{v}_{0}$ and   $z\left(\overline{v}_{0}\right)=0,$ we can easily calculate the value of derivative $z^{\prime}(\overline{v}_{0})$
\begin{eqnarray}\label{D_54}
    z^{\prime}\left(v_{0}\right)&=&\lim\limits_{f\to \overline{v}_{0}}\frac{z\left(f\right)-z\left(\overline{v}_{0}\right)}{f-\overline{v}_{0}}=\\
    &=&\lim\limits_{f\to \overline{v}_{0}}\frac{1-(f-\overline{v}_0)}{(1+\overline{v}_0)[\|q\|_{0,1}f+\overline{N}(f)+\overline{N}^{\prime}\left(\overline{v}_{0}\right)\overline{v}_{0}-\overline{N}\left(\overline{v}_{0}\right)]}=\nonumber\\
    &=&\frac{1}{\overline{v}_{0}(1+\overline{v}_0)[\|q\|_{0,1}+\overline{N}^{\prime}\left(\overline{v}_{0}\right)]}>0.
    \nonumber
\end{eqnarray}
Inequality \eqref{D_54} implies that there exists an inverse function $f=z^{-1}$ which is holomorphic in some interval $(-R,R)$ (see \cite[p. 87]{Analytic_theory}). Now let us prove that series \eqref{45} converges at the endpoint $z=R.$ Suppose that this is not the case and  series \eqref{45} diverges at the point $z=R,$ that is,
 $$\lim\limits_{z\to R-0}f(z)=+\infty.$$
 However, taking into account that equality \eqref{D_52} holds for every $z$ in $(-R, R)$, we immediately get the contradiction
 \begin{eqnarray}
   1 = \lim\limits_{z\to R-0}\left (f(z)-\overline{v}_{0}+\right . \nonumber\\
     \left. + \frac{(1+\overline{v}_0)z[\|q\|_{0,1}f(z)+\overline{N}(f(z))+\overline{N}^{\prime}\left(\overline{v}_{0}\right) \overline{v}_{0}-\overline{N}\left(\overline{v}_{0}\right)]}{f(z)-\overline{v}_{0}}\right ) = +\infty.
 \end{eqnarray}
   This contradiction implies that $f(R)<+\infty$. Thus, we have the inequality
$$R^{j}\overline{v}^j\leq\frac{c}{j^{1+\varepsilon}}$$
with some positive constants $c$ and $\varepsilon$, where the constant $R$ is determined by the value of $\|q\|_{0,1}$ and function $\overline{N}\left(u\right).$

Returning to notation (\ref{42}), we formally obtain estimates \eqref{u_estimate_1} and \eqref{lambda_estimate_1}.
Hence, for $n$ sufficiently large the numerical solution converges to the exact one with  error estimates \eqref{22}, where $R$ denotes the radius of convergence to the power series \eqref{D_52}, \eqref{45}.

Now we apply the same technique as above to the case when $\lambda_n^{(0)}=4\pi^2n^2$ $(n=1,2,\ldots )$.

If $n$ is an even number then
$$\lambda_n^{(j+1)}= \hspace*{100mm}$$
\begin{equation}\label{47}
=\frac{8\pi n}{3}
\left\{\int\limits_0^{1/2}\left [-\sum_{p=1}^j\lambda_n^{(j+1-p)}u_{n1}^{(p)}+q(x)u_{n1}^{(j)}(x)+ A_{j}(N; u_{n1}^{(0)},...,u_{n1}^{(j)})\right ]\sin\left(2\pi nx\right)dx-\right .
\end{equation}
$$\left .-\int\limits_{1/2}^1\left [-\sum_{p=1}^j\lambda_n^{(j+1-p)}u_{n2}^{(p)}+q(x)u_{n2}^{(j)}(x)+ A_{j}(N; u_{n2}^{(0)},...,u_{n2}^{(j)})\right ]{\sin\left(2\pi n x\right)}dx\right\}.$$
and functions $u_{ni}^{(j)}(x)$ $(i=1,2)$ can be found using formulas (\ref{13}), (\ref{37}) where parameter $c_2^{(j+1)}$ is determined by formulas (\ref{26}), (\ref{37}).

If $n$ is an odd number, then
$$\lambda_n^{(j+1)}= \hspace*{100mm}$$
\begin{equation}\label{48}
=8\pi n
\left\{\int\limits_0^{1/2}\left [-\sum_{p=1}^j\lambda_n^{(j+1-p)}u_{n1}^{(p)}+q(x)u_{n1}^{(j)}(x)+ A_{j}(N; u_{n1}^{(0)},...,u_{n1}^{(j)})\right ]\sin\left(2\pi nx\right)dx-\right .
\end{equation}
$$\left .-\int\limits_{1/2}^1\left [-\sum_{p=1}^j\lambda_n^{(j+1-p)}u_{n2}^{(p)}+q(x)u_{n2}^{(j)}(x)+ A_{j}(N; u_{n2}^{(0)},...,u_{n2}^{(j)})\right ]{\sin\left(2\pi n x\right)}dx\right\}.$$
and functions $u_{ni}^{(j)}(x)$ $(i=1,2)$ can be found using formulas  (\ref{13}), (\ref{37}) where parameter $c_2^{(j+1)}$ is determined by formulas (\ref{26}), (\ref{37}).

As a result we formally obtain estimates (\ref{28}) for an even number $n$, and estimates (\ref{31}) for an odd number $n,$ provided that
\begin{equation}\label{50}
r_n=\frac{1}{\pi nR}<1.
\end{equation}

Therefore, we arrive at the following convergence result for the nonlinear problem.

\begin{theorem}

\indent
\begin{enumerate}

\item[(a)]Suppose that the index $n$ of a trial eigenpair satisfies condition \eqref{23} where  $\lambda_n^{(0)}=4\pi^2(\pm\frac{2}{3}+2n)^2$, $(n=0,1,\ldots)$ and $R$ is a radius of convergence to series \eqref{45} satisfying equation (\ref{D_52}). Then FD-algorithm (\ref{7}), (\ref{11}), (\ref{37}), (\ref{40}) converges to the corresponding eigensolution of problem (\ref{33}) super-exponentially with error estimate (\ref{22}).
\item[(b)]Suppose that the index $n$ of a trial eigenpair satisfies condition (\ref{50}) where $R$ ia a radius of convergence to series \eqref{45} satisfying  equation (\ref{D_52}) and $\lambda_n^{(0)}=4\pi^2n^2$ $(n=1,2,\ldots)$.  Then the FD-algorithm (\ref{7}), (\ref{12}), (\ref{37}) and (\ref{47}) (if $n$ is even) or \eqref{48} (if $n$ is odd) converges to the corresponding eigensolution of problem (\ref{33}) super-exponentially with error estimate (\ref{28}) for even $n$ and with error estimate (\ref{31}) for odd $n$.
\end{enumerate}
\end{theorem}

\section{Numerical examples}\label{s_4}

\subsection{Example 1}

As an example we consider the following nonlinear eigenvalue problem

$$u_1''(x)-(x+3x^2)u_1(x)+\lambda u_1(x)-\left[u_1(x)\right]^{2}=0, \hspace*{5mm} x\in \left(0,\frac{1}{2}\right)$$
\begin{equation}\label{ex_1}
u_2''(x)-(x+3x^2)u_2(x)+\lambda u_2(x)-\left[u_2(x)\right]^{2}=0, \hspace*{5mm} x\in \left(\frac{1}{2},1\right)
\end{equation}
$$\left [u\left (\frac{1}{2}\right )\right ]=0, \hspace*{5mm} \left [u'\left (\frac{1}{2}\right )\right ]=1,$$
$$u_1(0)=u_2(1)=0, \hspace*{5mm} u'(0)=1.$$

Let us compute the approximations to the six least eigenvalues of the problem. From formulas \eqref{11} and \eqref{12} we obtain the zero approximations for the FD-method's algorithm:

from \eqref{11} with $n=0:$
\begin{equation}\label{ex_2}
    \lambda_{0}^{(0)}=\frac{16}{9}\pi^2,\quad u_{0,1}^{(0)}\left(x\right)=\frac{\sin(\frac{4}{3}\pi x)}{\frac{4}{3}\pi},\quad u_{0,2}^{(0)}\left(x\right)=\frac{\sin\left(\frac{4}{3}\pi\left(1-x\right)\right)}{\frac{4}{3}\pi};
\end{equation}

from \eqref{12} with $n=1:$
\begin{equation}\label{ex_3*}
    \lambda_{1}^{(0)}=4\pi^2,\quad u_{1,1}^{(0)}\left(x\right)=\frac{\sin(2\pi x)}{2\pi},\quad u_{1,2}^{(0)}\left(x\right)=0;
\end{equation}

from \eqref{11} with $n=1:$
\begin{equation}\label{ex_4*}
    \lambda_{2}^{(0)}=\frac{64}{9}\pi^2,\quad u_{2, 1}^{(0)}\left(x\right)=\frac{\sin(\frac{8}{3}\pi x)}{\frac{8}{3}\pi},\quad u_{2,2}^{(0)}\left(x\right)=\frac{\sin\left(\frac{8}{3}\pi\left(1-x\right)\right)}{\frac{8}{3}\pi};
\end{equation}

from \eqref{12} with $n=2:$
\begin{equation}\label{ex_5}
    \lambda_{3}^{(0)}=16\pi^2,\quad u_{3, 1}^{(0)}\left(x\right)=\frac{\sin(4\pi x)}{4\pi},\quad u_{3,2}^{(0)}\left(x\right)=2\frac{\sin\left(4\pi\left(x\right)\right)}{4\pi};
\end{equation}

from \eqref{11} with $n=1:$
\begin{equation}\label{ex_6}
    \lambda_{4}^{(0)}=\frac{256}{9}\pi^2,\quad u_{4, 1}^{(0)}\left(x\right)=\frac{\sin(\frac{16}{3}\pi x)}{\frac{16}{3}\pi},\quad u_{4,2}^{(0)}\left(x\right)=\frac{\sin\left(\frac{16}{3}\pi\left(1-x\right)\right)}{\frac{16}{3}\pi};
\end{equation}

for \eqref{12} with $n=3:$
\begin{equation}\label{ex_6*}
    \lambda_{5}^{(0)}=36\pi^2,\quad u_{5, 1}^{(0)}\left(x\right)=\frac{\sin(6\pi x)}{6\pi},\quad u_{5,2}^{(0)}\left(x\right)=0.
\end{equation}

We have calculated the FD-approximations for the eigenvalues $\lambda_n$ and corresponding eigenfunctions $u_{n}(x),$  $n=\overline{0,5}$ up to the rank $m=4$ inclusive. Since the exact solution of problem \eqref{ex_1} is unknown, in order to estimate the error, we analyzed the $L_\infty$-norm of the residual
\begin{equation}\label{ex_4}
    \overset{m}{\nu}_{ni}(x)=\frac{d^{2}}{dx^{2}}\overset{m}{u}_{ni}(x)-(x+3x^2)\overset{m}{u}_{ni}(x)+\overset{m}{\lambda}_{n} \overset{m}{u}_{ni}(x)-\left[\overset{m}{u}_{ni}(x)\right]^{2}.
\end{equation}

\begin{figure}[htbp]
\begin{minipage}[h]{1\linewidth}
\begin{minipage}[h]{0.48\linewidth}
\center{\rotatebox{0}{\includegraphics[width=1\linewidth]{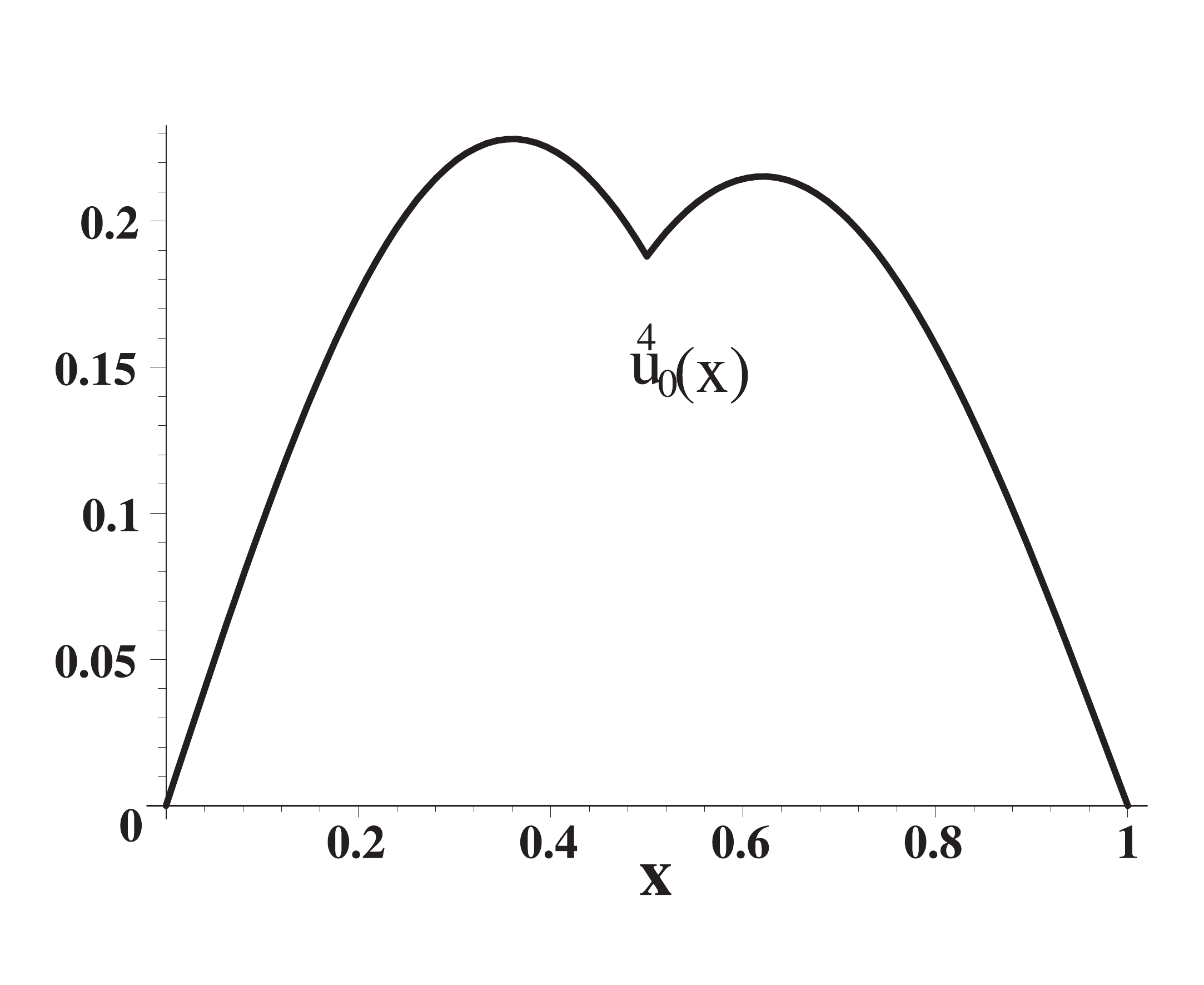}} \\ a) }
\end{minipage}
\hfill\label{image1}
\begin{minipage}[h]{0.48\linewidth}
\center{\rotatebox{0}{\includegraphics[width=1\linewidth]{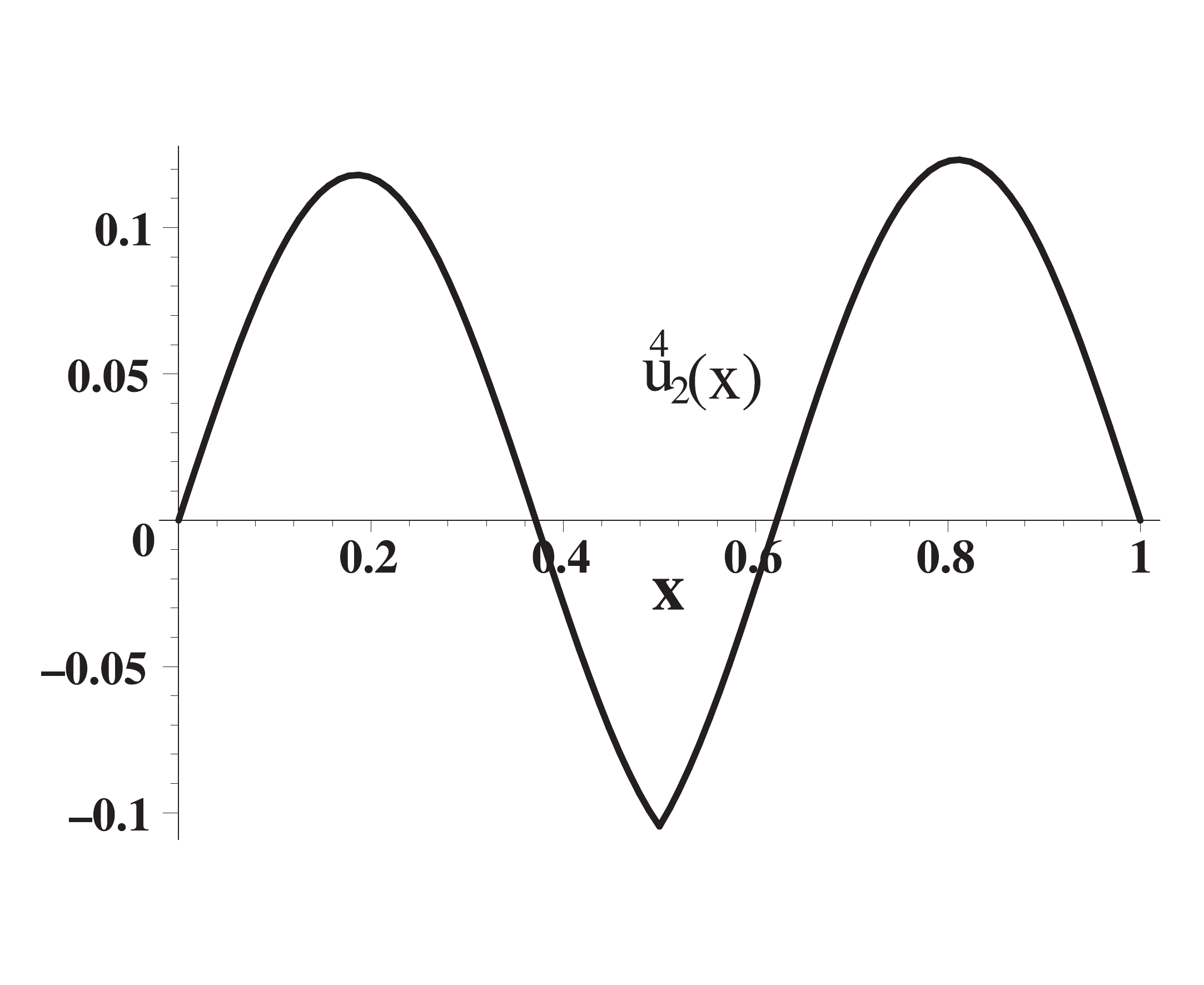}} \\ b) }
\end{minipage}
\caption{Example 4.1. FD-approximation of $u_{0}(x)$ (a) and $u_{2}(x)$ (b).}\label{fig_1}
\end{minipage}
\end{figure}

\begin{figure}[htbp]
\begin{minipage}[h]{1\linewidth}
\begin{minipage}[h]{0.48\linewidth}
\center{\rotatebox{-0}{\includegraphics[width=1\linewidth]{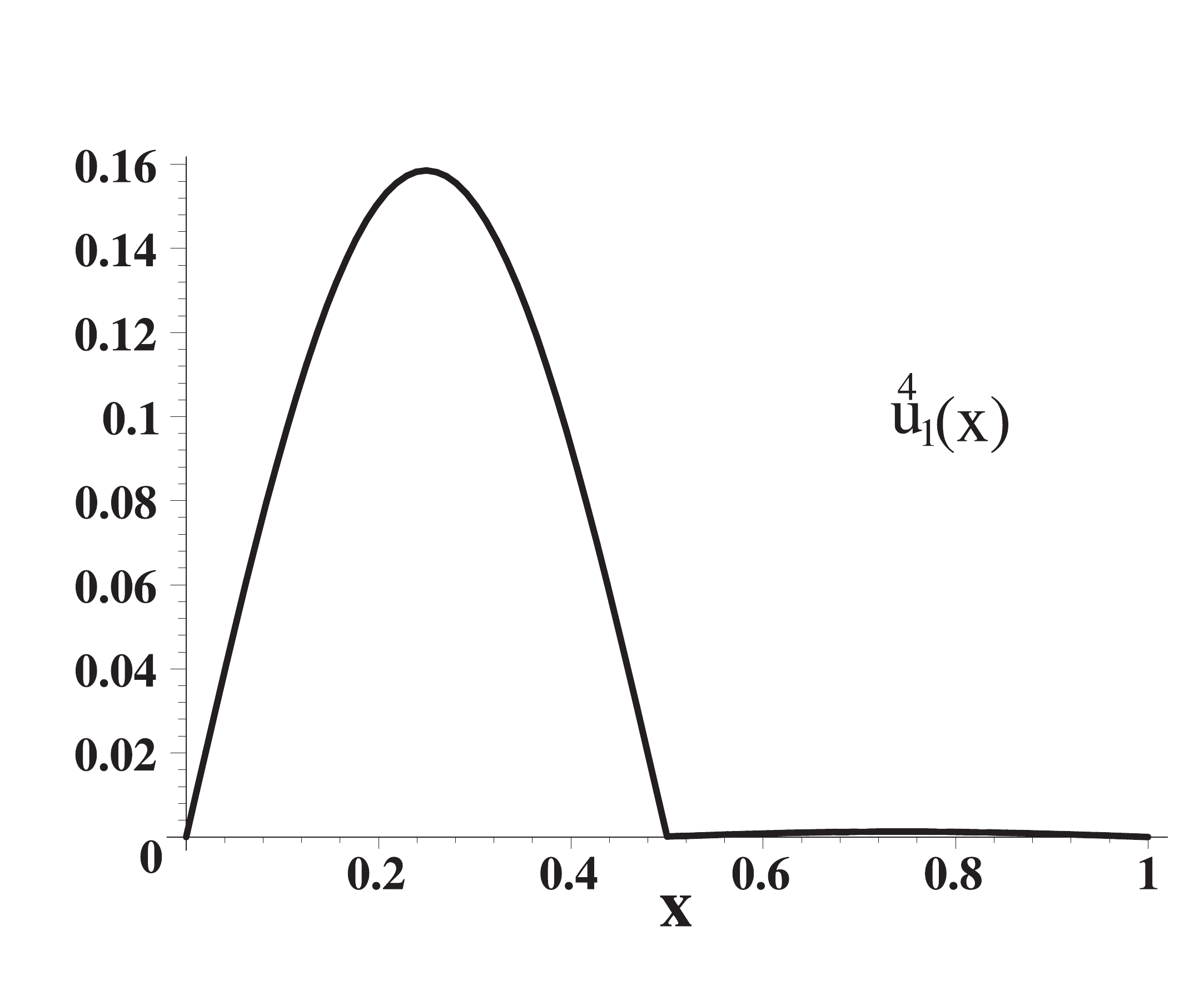}} \\ a) }
\end{minipage}
\hfill\label{image2}
\begin{minipage}[h]{0.48\linewidth}
\center{\rotatebox{-0}{\includegraphics[width=1\linewidth]{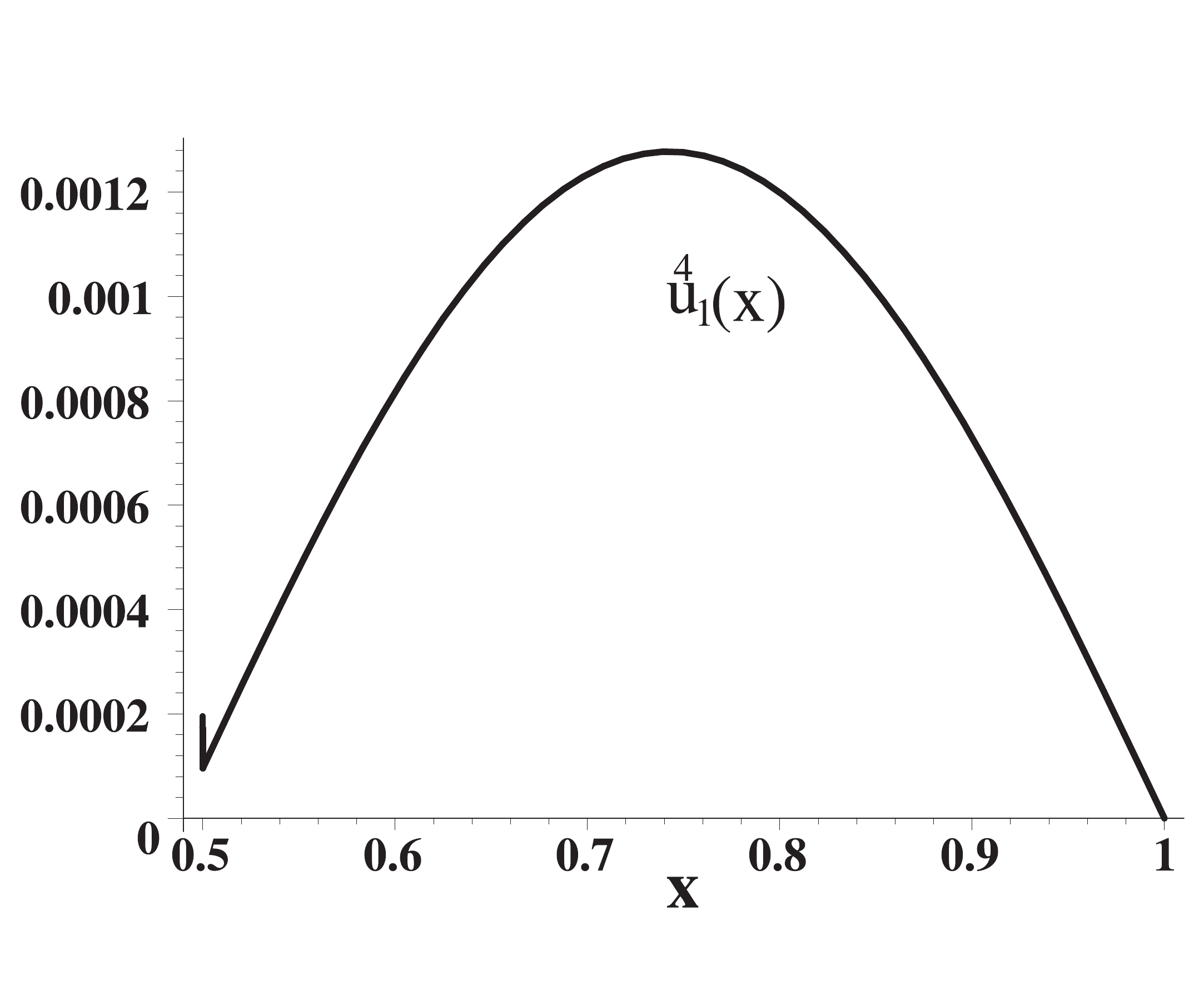}} \\ b) }
\end{minipage}
\caption{Example 4.1 FD-approximation of $u_{1}(x)$ in $[0,1]$ (a) and in $(0.5,1]$ (b)}\label{fig_2}
\end{minipage}
\end{figure}

\begin{figure}[htbp]
\begin{minipage}[h]{1\linewidth}
\begin{minipage}[h]{0.48\linewidth}
\center{\rotatebox{-0}{\includegraphics[width=1\linewidth]{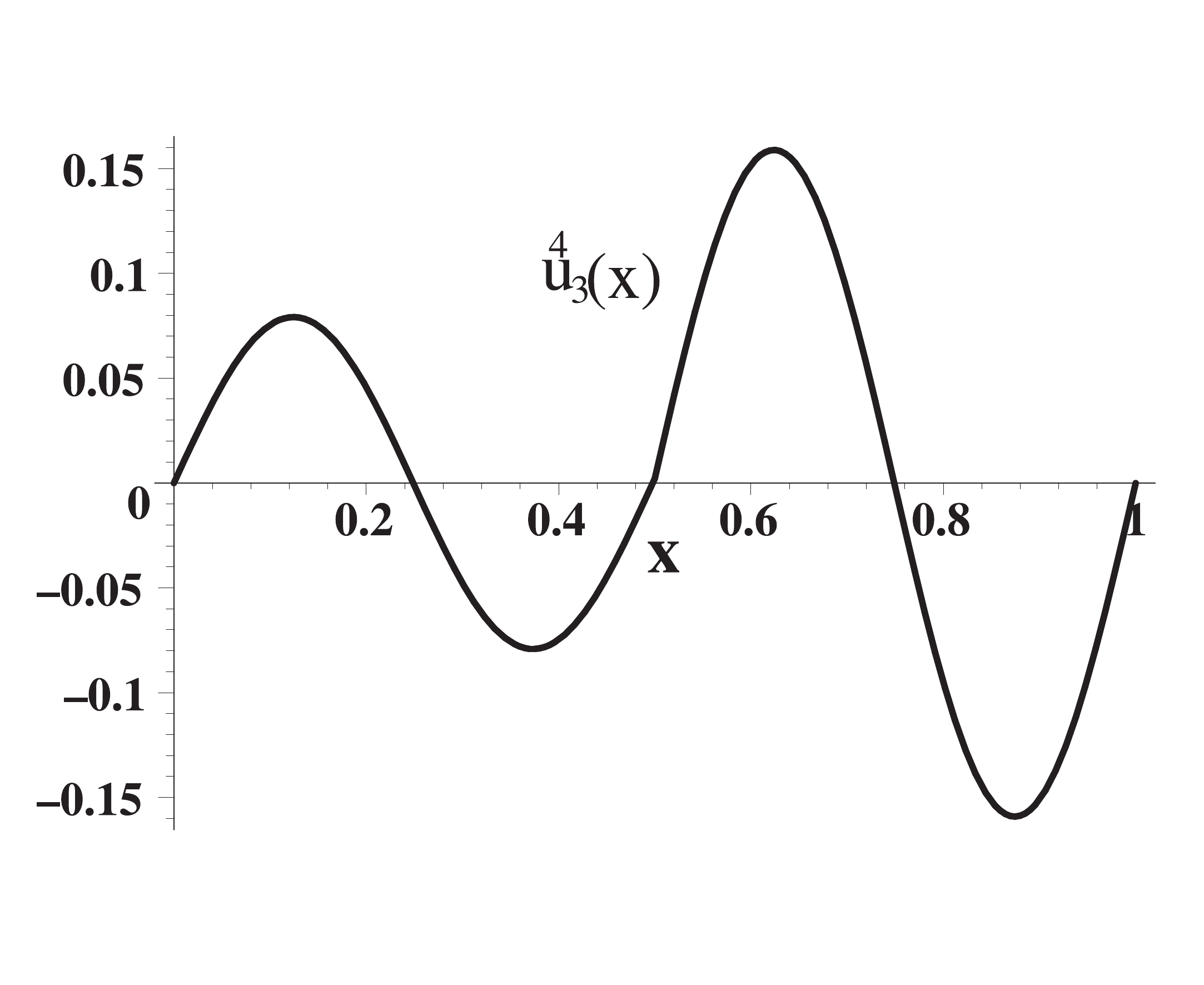}} \\  a)}
\end{minipage}
\hfill
\begin{minipage}[h]{0.48\linewidth}
\center{\rotatebox{-0}{\includegraphics[width=1\linewidth]{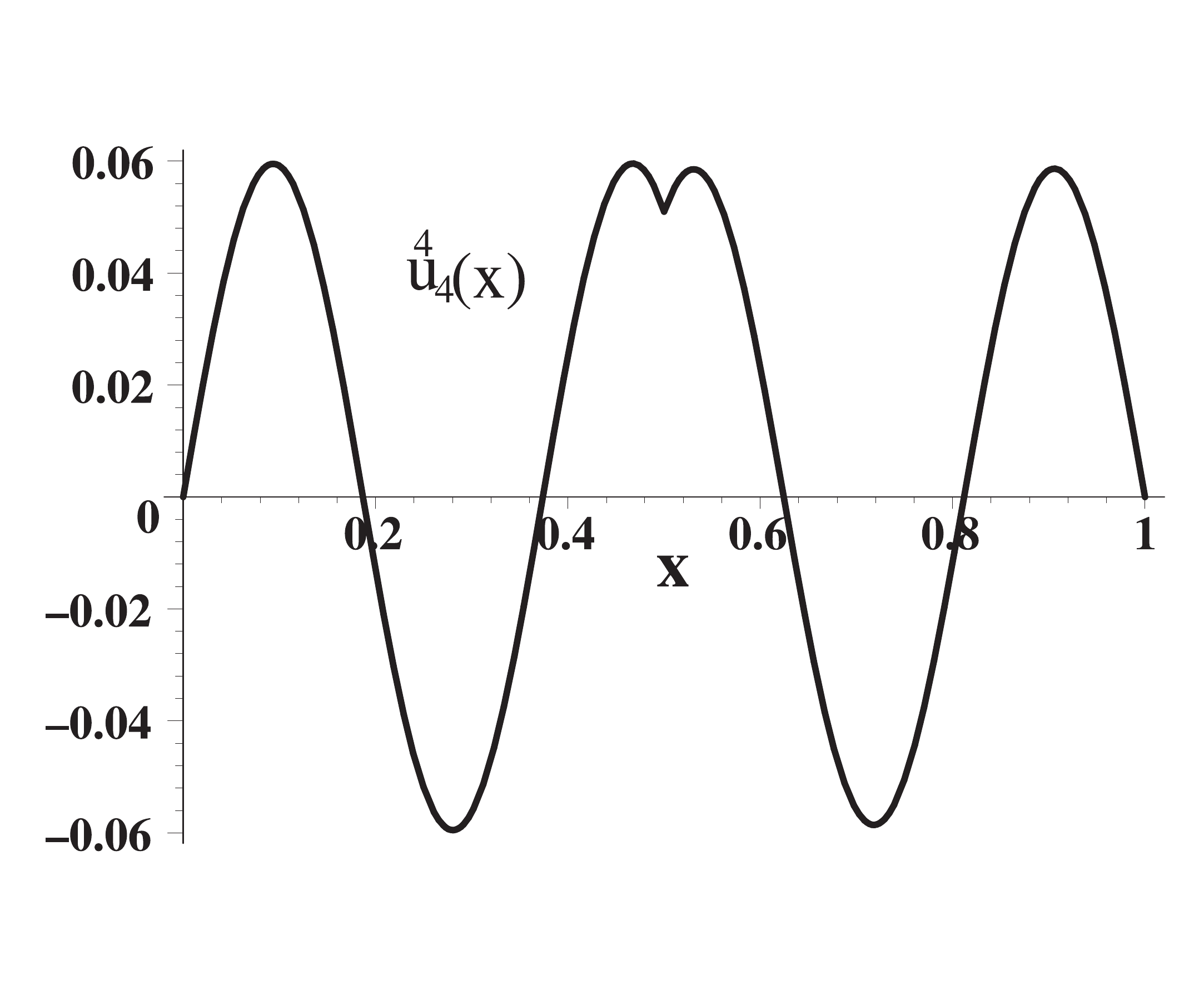}} \\ b)}
\end{minipage}
\caption{Example 4.1. FD-approximation of $u_{3}(x)$ (a)) and $u_{4}(x)$ (b)).}\label{fig_3}
\end{minipage}
\end{figure}

\begin{figure}[htbp]
\begin{minipage}[h]{1\linewidth}
\begin{minipage}[h]{0.48\linewidth}
\center{\rotatebox{-0}{\includegraphics[width=1\linewidth]{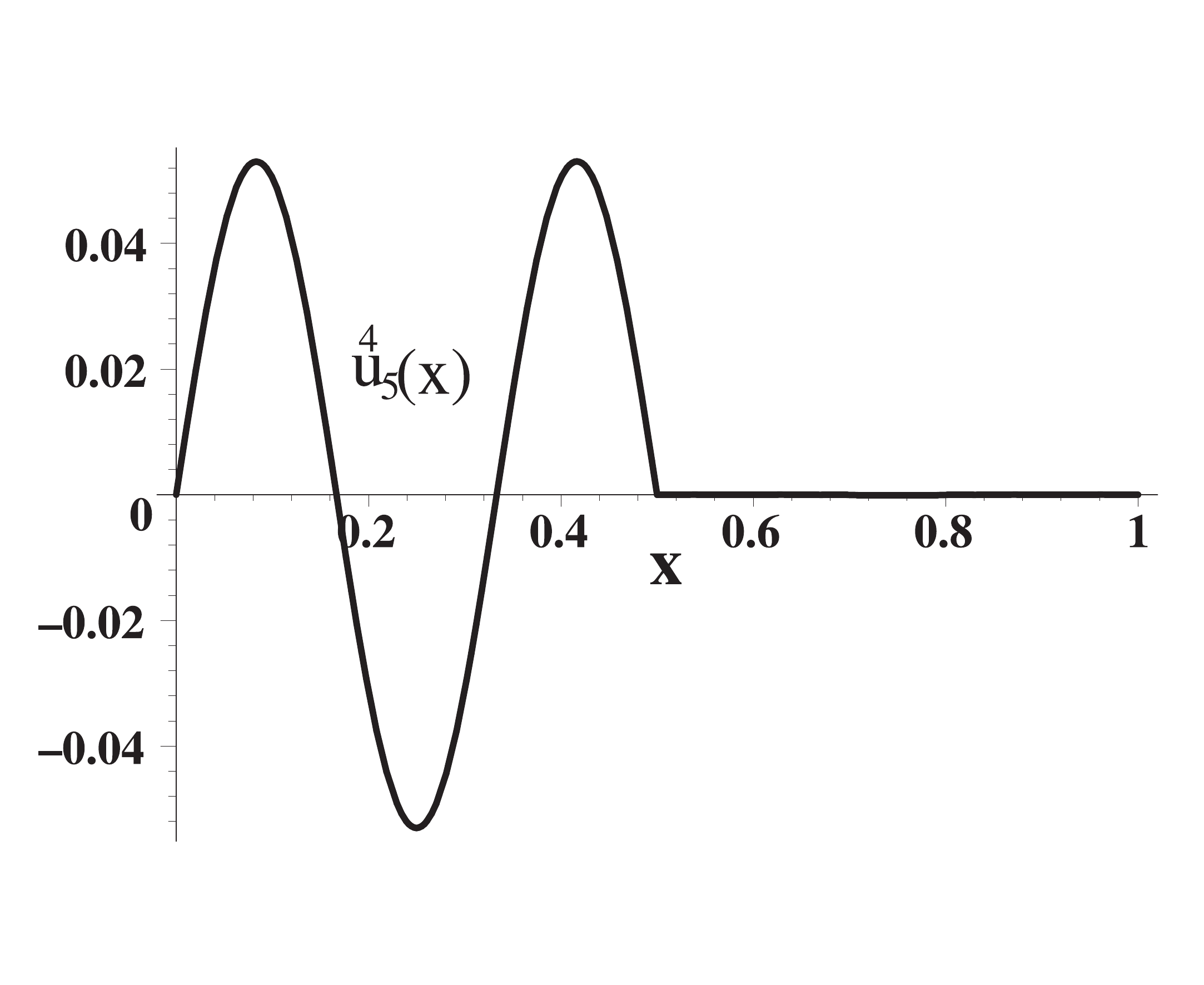}} \\ a) }
\end{minipage}
\hfill\label{image4}
\begin{minipage}[h]{0.48\linewidth}
\center{\rotatebox{-0}{\includegraphics[width=1.0\linewidth]{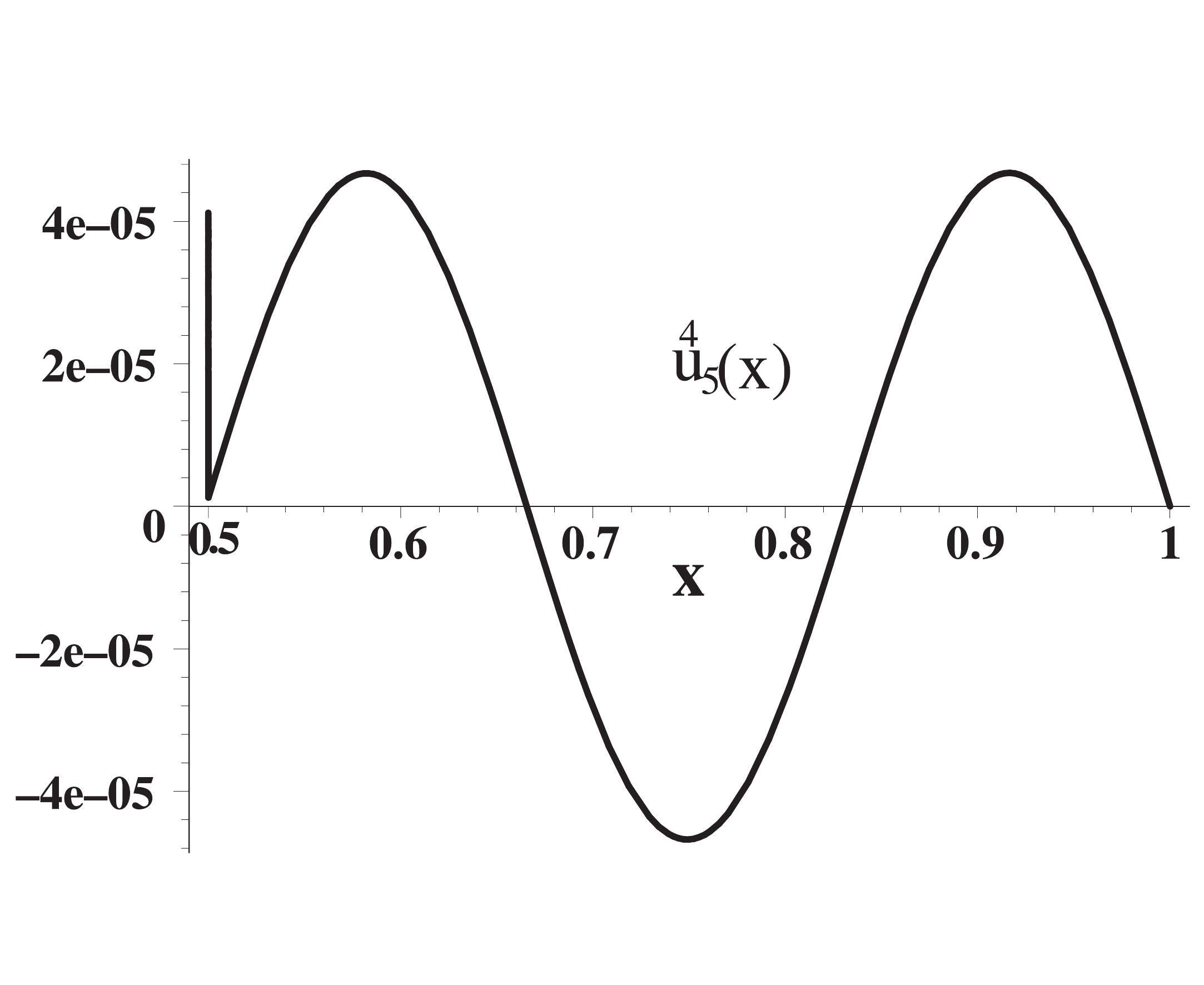}} \\ b) }
\end{minipage}
\caption{Example 4.1. FD-approximation of $u_{5}(x)$ in $[0,1]$ (a) and in $(0.5,1]$ (b).}\label{fig_4}
\end{minipage}
\end{figure}

\begin{figure}[htbp]
\begin{minipage}[h]{1\linewidth}
\center{\rotatebox{-0}{\includegraphics[width=0.5\linewidth]{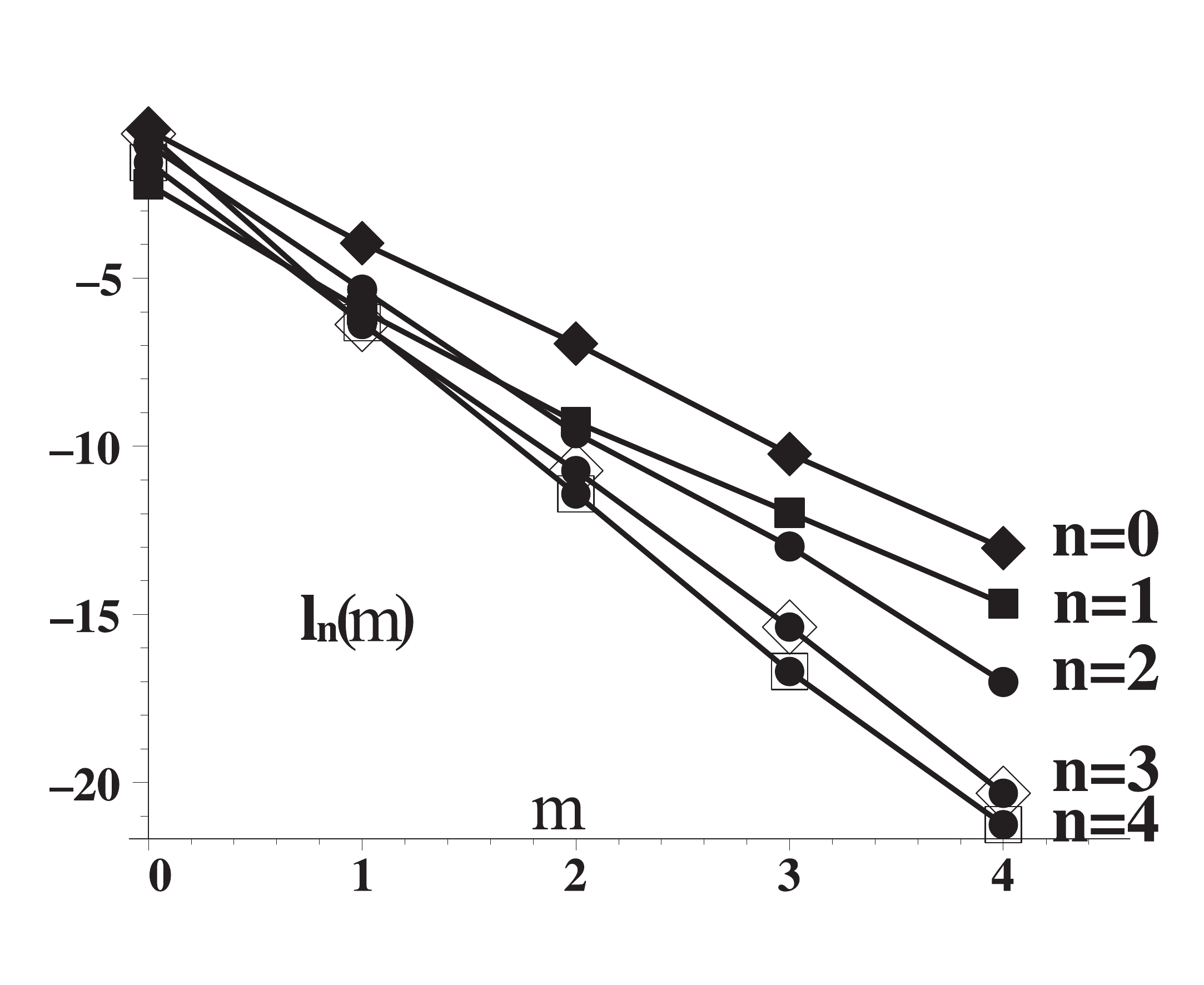}} \\ }
\end{minipage}
\caption{Example 4.1. $L_\infty$-norm of residual \eqref{ex_4} versus the rank $m$ of the approximation on a semi-logarithmic scale.}\label{fig_5}
\end{figure}

The results obtained are presented in tables \ref{tabl_1} -- \ref{tabl_9} and depicted in figures \ref{fig_1} -- \ref{fig_5}.

\begin{table}[htbp]
\caption{Example 4.1. FD-approximation of $\lambda_{0}$ and $u_{0}(x).$}\label{tabl_1}
\begin{tabular}{|c|c|c|c|c|}
  \hline
  $m$ & $\overset{m}{\lambda}_{0}$ & $\|u_{0, 1}^{(m)}(x)\|_{\infty, [0,1/2]}$ & $\|u_{0, 2}^{(m)}(x)\|_{\infty, [1/2, 1]}$ & $\|\nu^{(m)}_{0}(x)\|_{\infty, [0,1]}$ \\
  \hline
  0 & 17.5459633797144 & 0.24 & 0.24 & 0.55 \\
  \hline
  1 & 19.6940699073641 & 0.19e-1 & 0.25e-1 & 0.19e-1 \\
  \hline
  2 & 19.6740808021547 & 0.60e-3 & 0.11e-2 & 0.96e-3 \\
  \hline
  3 & 19.6754846046439 & 0.28e-4 & 0.51e-4 & 0.36e-4 \\
  \hline
  4 & 19.6754786167117 & 0.88e-6 & 0.20e-5 & 0.22e-5 \\
  \hline
\end{tabular}
\newline
$\overset{4}{\lambda_{0}}= 19.6754786167117.$
\end{table}

\begin{table}[htbp]
\caption{Example 4.1. FD-approximation of $\lambda_{1}$ and $u_{1}(x).$}\label{tabl_2}
\begin{tabular}{|c|c|c|c|c|}
  \hline
  $m$ & $\overset{m}{\lambda}_{1}$ & $\|u_{1, 1}^{(m)}(x)\|_{\infty, [0,1/2]}$ & $\|u_{1, 2}^{(m)}(x)\|_{\infty, [1/2, 1]}$ & $\|\nu^{(m)}_{1}(x)\|_{\infty, [0,1]}$ \\
  \hline
  0 & 39.4784176043574 & 0.16 & 0. & 0.11 \\
  \hline
  1 & 40.0755170720146 & 0.80e-3 & 0.13e-2 & 0.26e-2 \\
  \hline
  2 & 40.0597952320099 & 0.94e-4 & 0.93e-4 & 0.93e-4 \\
  \hline
  3 & 40.0595106843757 & 0.20e-5 & 0.29e-5 & 0.63e-5 \\
  \hline
  4 & 40.0594734299829 & 0.27e-6 & 0.27e-6 & 0.42e-6 \\
  \hline
\end{tabular}
\newline
$\overset{4}{\lambda_{1}}=40.0594734299829.$
\end{table}

\begin{table}[htbp]
\caption{Example 4.1. FD-approximation of $\lambda_{2}$ and $u_{2}(x).$}\label{tabl_3}
\begin{tabular}{|c|c|c|c|c|}
  \hline
  $m$ & $\overset{m}{\lambda}_{2}$ & $\|u_{2, 1}^{(m)}(x)\|_{\infty, [0,1/2]}$ & $\|u_{2, 2}^{(m)}(x)\|_{\infty, [1/2, 1]}$ & $\|\nu^{(m)}_{2}(x)\|_{\infty, [0,1]}$ \\
  \hline
  0 & 70.1838535188575 & 0.12 & 0.12 & 0.36 \\
  \hline
  1 & 71.9666624409054 & 0.37e-2 & 0.46e-2 & 0.48e-2 \\
  \hline
  2 & 71.9768253657124 & 0.53e-5 & 0.86e-4 & 0.67e-4 \\
  \hline
  3 & 71.9766663564735 & 0.12e-5 & 0.20e-5 & 0.23e-5 \\
  \hline
  4 & 71.9766690902938 & 0.42e-7 & 0.46e-7 & 0.41e-7 \\
  \hline
\end{tabular}
\newline
$\overset{4}{\lambda_{2}}= 71.9766690902938.$
\end{table}

\begin{table}[htbp]
\caption{Example 4.1. FD-approximation of $\lambda_{3}$ and $u_{3}(x).$}\label{tabl_4}
\begin{tabular}{|c|c|c|c|c|}
  \hline
  $m$ & $\overset{m}{\lambda}_{3}$ & $\|u_{3, 1}^{(m)}(x)\|_{\infty, [0,1/2]}$ & $\|u_{3, 2}^{(m)}(x)\|_{\infty, [1/2, 1]}$ & $\|\nu^{(m)}_{3}(x)\|_{\infty, [0,1]}$ \\
  \hline
  0 & 157.913670417429 & 0.80e-1 & 0.16 & 0.49 \\
  \hline
  1 & 159.737504889796 & 0.13e-2 & 0.20e-2 & 0.17e-2 \\
  \hline
  2 & 159.739358889030 & 0.18e-4 & 0.18e-4 & 0.22e-4 \\
  \hline
  3 & 159.739350000888 & 0.22e-6 & 0.21e-6 & 0.21e-6 \\
  \hline
  4 & 159.739350058922 & 0.48e-9 & 0.85e-9 & 0.15e-8 \\
  \hline
\end{tabular}
\newline
$\overset{4}{\lambda_{3}}= 159.739350058922.$
\end{table}

\begin{table}[htbp]
\caption{Example 4.1. FD-approximation of $\lambda_{4}$ and $u_{4}(x).$}\label{tabl_5}
\begin{tabular}{|c|c|c|c|c|}
  \hline
  $m$ & $\overset{m}{\lambda}_{4}$ & $\|u_{4, 1}^{(m)}(x)\|_{\infty, [0,1/2]}$ & $\|u_{4, 2}^{(m)}(x)\|_{\infty, [1/2, 1]}$ & $\|\nu^{(m)}_{4}(x)\|_{\infty, [0,1]}$ \\
  \hline
  0 & 280.735414075430 & 0.60e-1 & 0.60e-1 & 0.21 \\
  \hline
  1 & 282.620725160874 & 0.57e-3 & 0.13e-2 & 0.18e-2 \\
  \hline
  2 & 282.622515077477 & 0.12e-4 & 0.10e-4 & 0.11e-4 \\
  \hline
  3 & 282.622528329338 & 0.36e-7 & 0.79e-7 & 0.56e-7 \\
  \hline
  4 & 282.622528620046 & 0.28e-9 & 0.24e-9 & 0.59e-9 \\
  \hline
\end{tabular}
\newline
$\overset{4}{\lambda_{4}}= 282.622528620046.$
\end{table}

\begin{table}[htbp]
\caption{Example 4.1 FD-approximation of $\lambda_{5}$ and $u_{5}(x).$}\label{tabl_6}
\begin{tabular}{|c|c|c|c|c|}
  \hline
  $m$ & $\overset{m}{\lambda}_{5}$ & $\|u_{5, 1}^{(m)}(x)\|_{\infty, [0,1/2]}$ & $\|u_{5, 2}^{(m)}(x)\|_{\infty, [1/2, 1]}$ & $\|\nu^{(m)}_{5}(x)\|_{\infty, [0,1]}$ \\
  \hline
  0 & 355.305758439216 & 0.53e-1 & 0. & 0.53e-1 \\
  \hline
  1 & 355.816547268956 & 0.11e-3 & 0.47e-4 & 0.14e-3 \\
  \hline
  2 & 355.814884402534 & 0.83e-6 & 0.10e-5 & 0.15e-5 \\
  \hline
  3 & 355.814878976396 & 0.24e-8 & 0.96e-8 & 0.47e-7 \\
  \hline
  4 & 355.814878544097 & 0.20e-9 & 0.30e-9 & 0.59e-9 \\
  \hline
\end{tabular}
\newline
$\overset{4}{\lambda_{5}}= 355.814878544097.$
\end{table}

Analyzing the results of calculations we have noticed that there exist eigenfunctions with some indexes $n$ such that their numerical approximation $\overset{4}{u}_n(x)$ have exactly $n-1$ zeros in the interval $(0,1)$ (see figure \ref{fig_2} for $n=1$ and  figure \ref{fig_3} for $n=5$). Zero approximations of such eigenpairs are determined by formula \eqref{12} with odd $n$. Thus, taking into account the numerical calculations and convergence result, we can conclude that there does not exist an eigenfunction of problem \eqref{ex_1} with exactly one or five zeros in the interval $(0,1)$. Further numerical experiments suggests that  there does not exist an eigenfunction which has exactly $1+4k$ $(k=0,1,2,\ldots)$ zeros in the interval $(0,1)$. Therefore, there exist two linearly independent eigenfunctions with exactly $4k$ zeros.

Hence, we can conclude that the theoretical approach developed in \cite{heinz4} and \cite{zhidkov} for nonlinear eigenvalue problems is not applicable for the problems of type \eqref{3}. The reason is that the corresponding differential operator is nonself-adjoint. So we see the crucial importance of the requirement that the differential operator is self-adjoint for the applicability of Theorem 4.2 and Theorem 5.6 from \cite{heinz4}.

To illustrate the exponential convergence rate of the method it is convenient to consider the function
\begin{equation}\label{log_funct}
   l=l_n(m)=\ln(\|\overset{m}{\nu}_{n}\|_{\infty}).
\end{equation}
  The values presented in table \ref{tabl_9} suggests that the graphs of function \eqref{log_funct} for $n\in\overline{0,5
  }$ are very similar to the straight lines, see  Figure \ref{fig_5}. And it is this fact that confirms the exponential nature of the FD-method's convergence. Also, it is easy to see that the slope of the lines on Figure \ref{fig_5}  increases as the index $n$ of the trial eigenvalue increases. This means that the convergence rate of the method  increases  along with the index of eigenvalue.

\begin{table}[htbp]
\caption{Example 4.1. The values of function $l_n(m)$ \eqref{log_funct}}\label{tabl_9}
\begin{tabular}{|c|c|c|c|c|c|c|}
  \hline
   & $n=0$ & $n=1$ & $n=2$ & $n=3$ & $n=4$ & $n=5$  \\
  \hline
  $m=0$ & -0.60 & -2.21 & -1.02 & -0.71 & -1.56 & -2.94\\
  \hline
  $m=1$ & -3.96 & -5.95 & -5.34 & -6.38 & -6.32 & -8.87\\
  \hline
  $m=2$ & -6.95 & -9.28 & -9.61 & -10.72 & -11.42 & -13.41\\
  \hline
  $m=3$ & -10.23 & -11.97 & -12.98 & -15.38 & -16.70 & -16.87\\
  \hline
  $m=4$ & -13.03 & -14.68 & -17.01 & -20.32 & -21.25 & -21.25\\
  \hline
\end{tabular}
\end{table}

\subsection{Example 2}

As the second example let us consider the following nonlinear eigenvalue transmission problem with singularity
$$u_1''(x)+\left (\lambda -(\frac{1}{2}-x)^{-1/2}\right )u_1(x)-\left[u_1(x)\right]^{2}=0, \hspace*{5mm} x\in \left(0,\frac{1}{2}\right),$$
\begin{equation}\label{ex_2_1}
u_2''(x)+\left (\lambda -(x-\frac{1}{2})^{-1/2}\right )u_2(x)-\left[u_2(x)\right]^{2}=0, \hspace*{5mm} x\in \left(\frac{1}{2},1\right),
\end{equation}
$$u_1(0)=u_2(1)=0, \hspace*{5mm} u_1'(0)=0.$$

It is easy to ensure that for the case when $q(x)=\cfrac{1}{\sqrt{\left|\frac{1}{2}-x\right|}}$ the indefinite integrals in formulas \eqref{13} cannot be expressed through the elementary functions. Therefore, to approximate functions $u_{ni}^{(k)}(x)$ we have used the numerical integration according to the Simpson's rule with the precision  $10^{-15}$ (for $n=0,1,2$) and $10^{-19}$ (for $n=3$)\footnote{The precision was estimated using a posteriori error analysis.}. Since the exact solution to the problem is unknown and $\overset{m}{u_{ni}}''(\frac{1}{2})=\infty$, for the error control we have used the norm of residual
$$\|\nu^{(m)}_{n}\|_{\infty(0,1)}=\max\left\{\|\nu^{(m)}_{n,1}\|_{\infty(0,\frac{1}{2})}, \|\nu^{(m)}_{n,2}\|_{\infty(\frac{1}{2},1)}\right\}$$
where
\begin{equation}\label{second_residual}
    \nu_{n, 1}^{(m)}(x)=\int\limits_{0}^{x}\left[\frac{d^{2}}{dx^{2}}\overset{m}{u}_{n1}(x)+ \left (\overset{m}{\lambda}_n-\frac{1}{\sqrt{\frac{1}{2}-x}}\right )\overset{m}{u}_{n1}(x)-(\overset{m}{u}_{n1}(x))^2\right] dx=
\end{equation}
$$=\frac{d}{d x}\overset{m}{u}_{n 1}(x)+\int\limits_{0}^{x}\left[\left (\overset{m}{\lambda}_n-\frac{1}{\sqrt{\frac{1}{2}-x}}\right )\overset{m}{u}_{n1}(x)-(\overset{m}{u}_{n1}(x))^2\right] dx, \hspace*{3mm} x\in (0,\frac{1}{2})$$

$$\nu_{n, 2}^{(m)}(x)=\nu_{n, 1}^{(m)}\left(\frac{1}{2}\right)+\int\limits_{\frac{1}{2}}^{x}\left[\frac{d^{2}}{dx^{2}}\overset{m}{u}_{n 2}(x)+ \left (\overset{m}{\lambda}_n-\frac{1}{\sqrt{x-\frac{1}{2}}}\right )\overset{m}{u}_{n 2}(x)-(\overset{m}{u}_{n 2}(x))^2\right] dx=$$
$$=\nu_{n, 1}^{(m)}\left(\frac{1}{2}\right)+\frac{d}{d x}\overset{m}{u}_{n 2}(x)-\frac{d}{d x}\overset{m}{u}_{n 2}(\frac{1}{2})+\int\limits_{\frac{1}{2}}^{x}\left[\left (\overset{m}{\lambda}_n-\frac{1}{\sqrt{x-\frac{1}{2}}}\right )\overset{m}{u}_{n 2}(x)-(\overset{m}{u}_{n 2}(x))^2\right] dx,$$
$x\in (\frac{1}{2},1)$.\\[2mm]

\begin{figure}[htbp]
\begin{minipage}[h]{1\linewidth}
\begin{minipage}[h]{0.48\linewidth}
\center{\rotatebox{-0}{\includegraphics[width=1\linewidth]{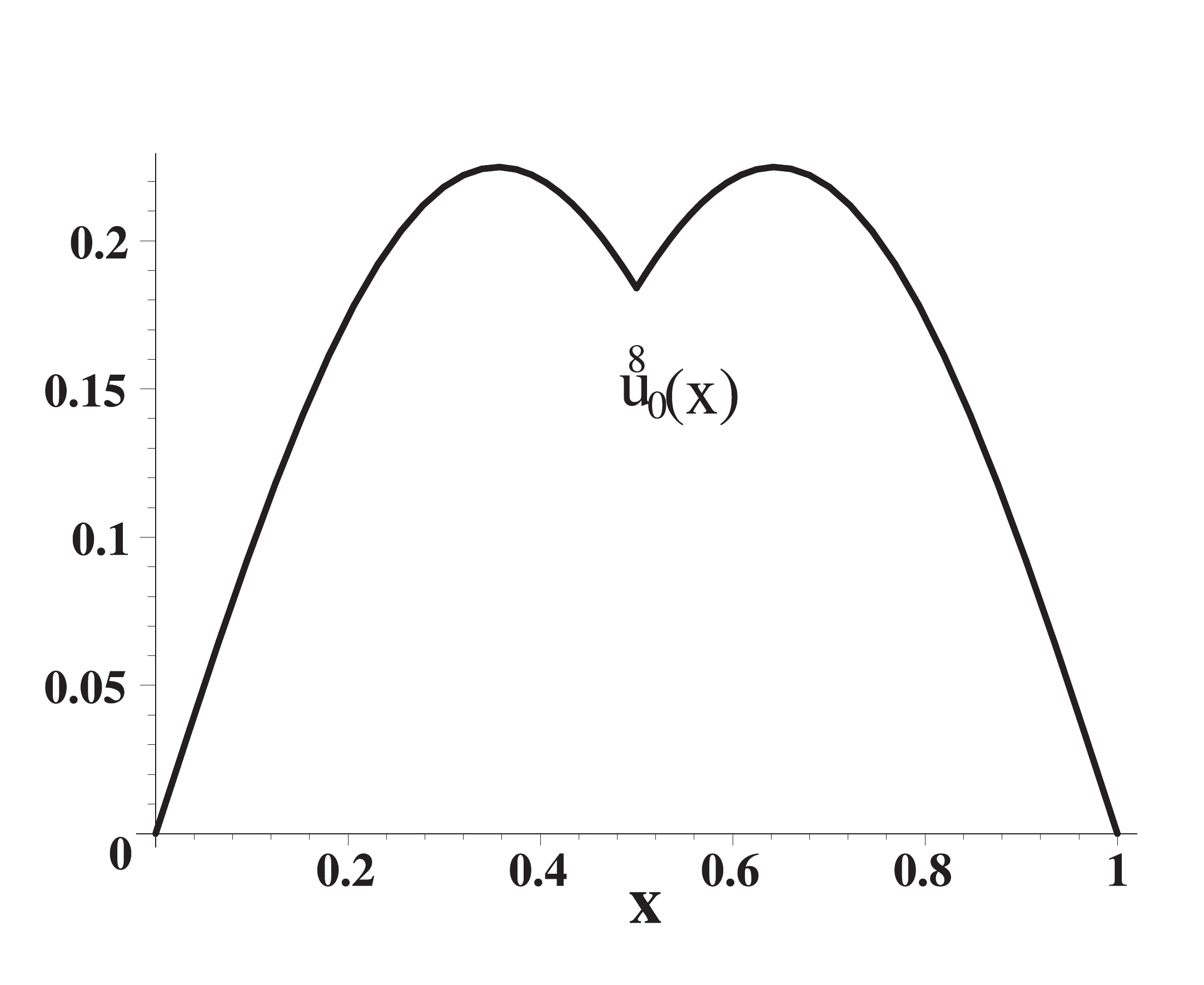}} \\ a) }
\end{minipage}
\hfill
\begin{minipage}[h]{0.48\linewidth}
\center{\rotatebox{-0}{\includegraphics[width=1\linewidth]{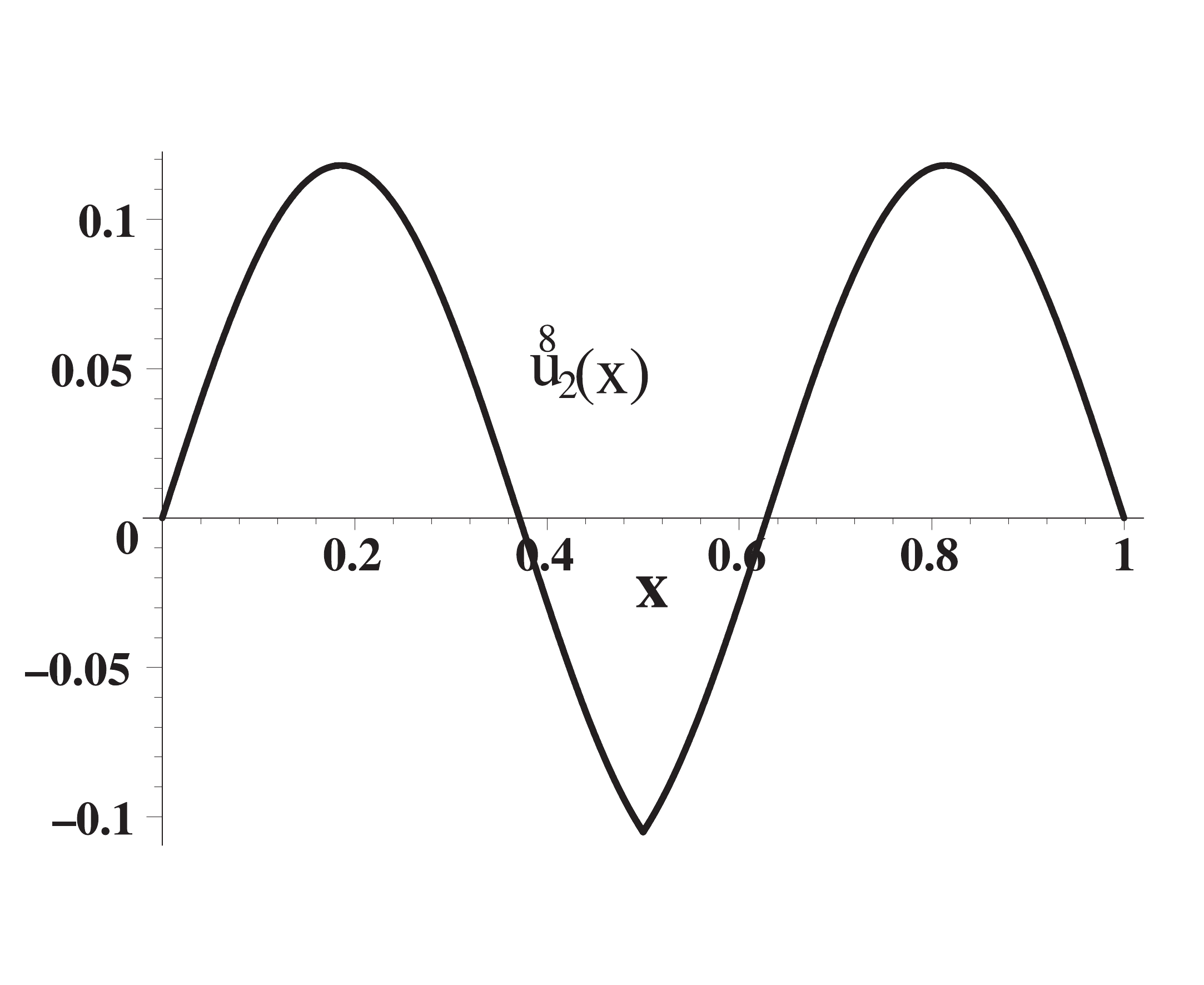}} \\ b) }
\end{minipage}
\caption{Example 4.2. FD-approximation of $u_{0}(x)$ (a) and $u_{2}(x)$ (b).} \label{fig_6}
\end{minipage}
\end{figure}

\begin{figure}[htbp]
\begin{minipage}[h]{1\linewidth}
\begin{minipage}[h]{0.48\linewidth}
\center{\rotatebox{-0}{\includegraphics[width=1\linewidth]{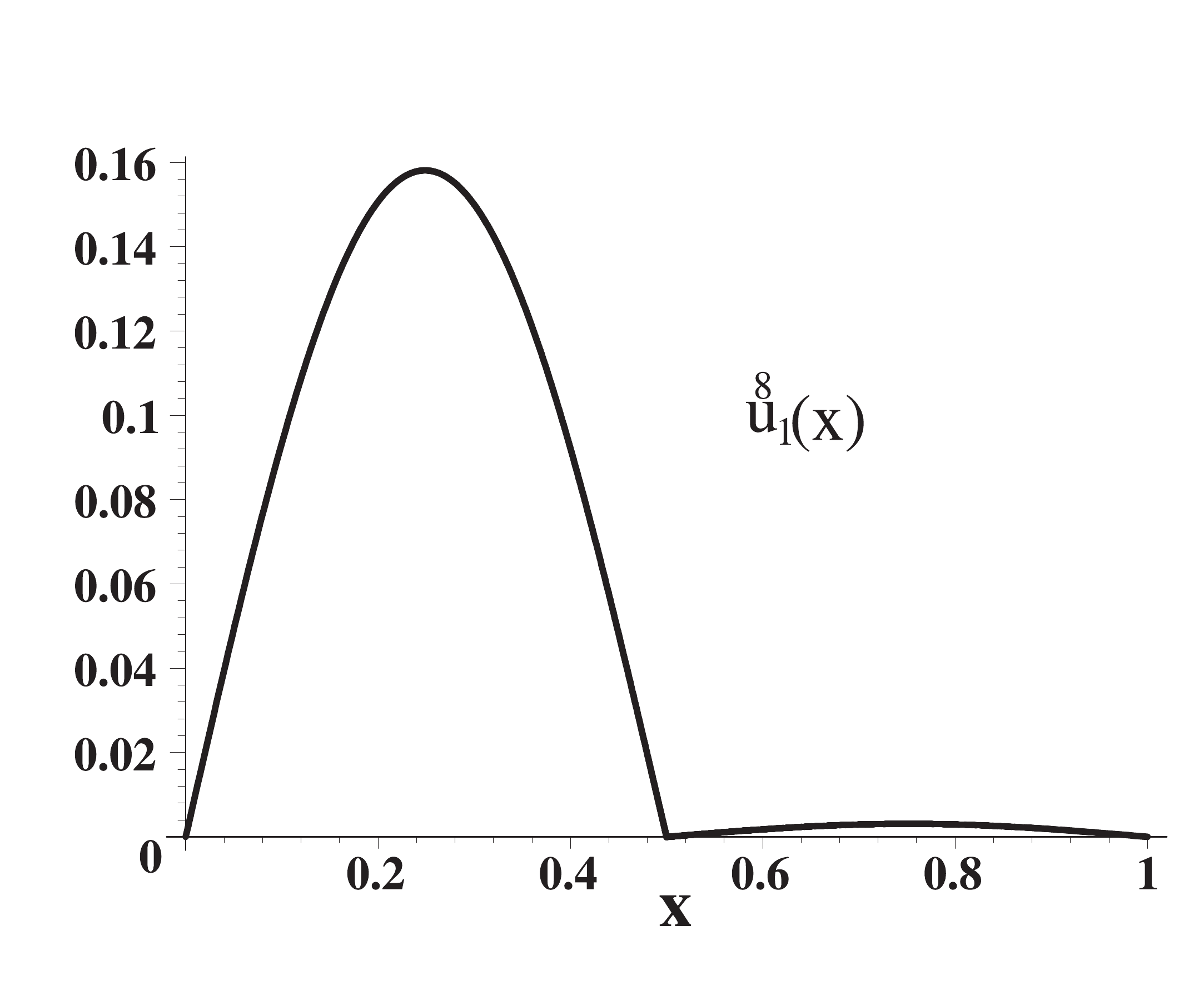}} \\ a) }
\end{minipage}
\hfill\label{image1}
\begin{minipage}[h]{0.48\linewidth}
\center{\rotatebox{-0}{\includegraphics[width=1\linewidth]{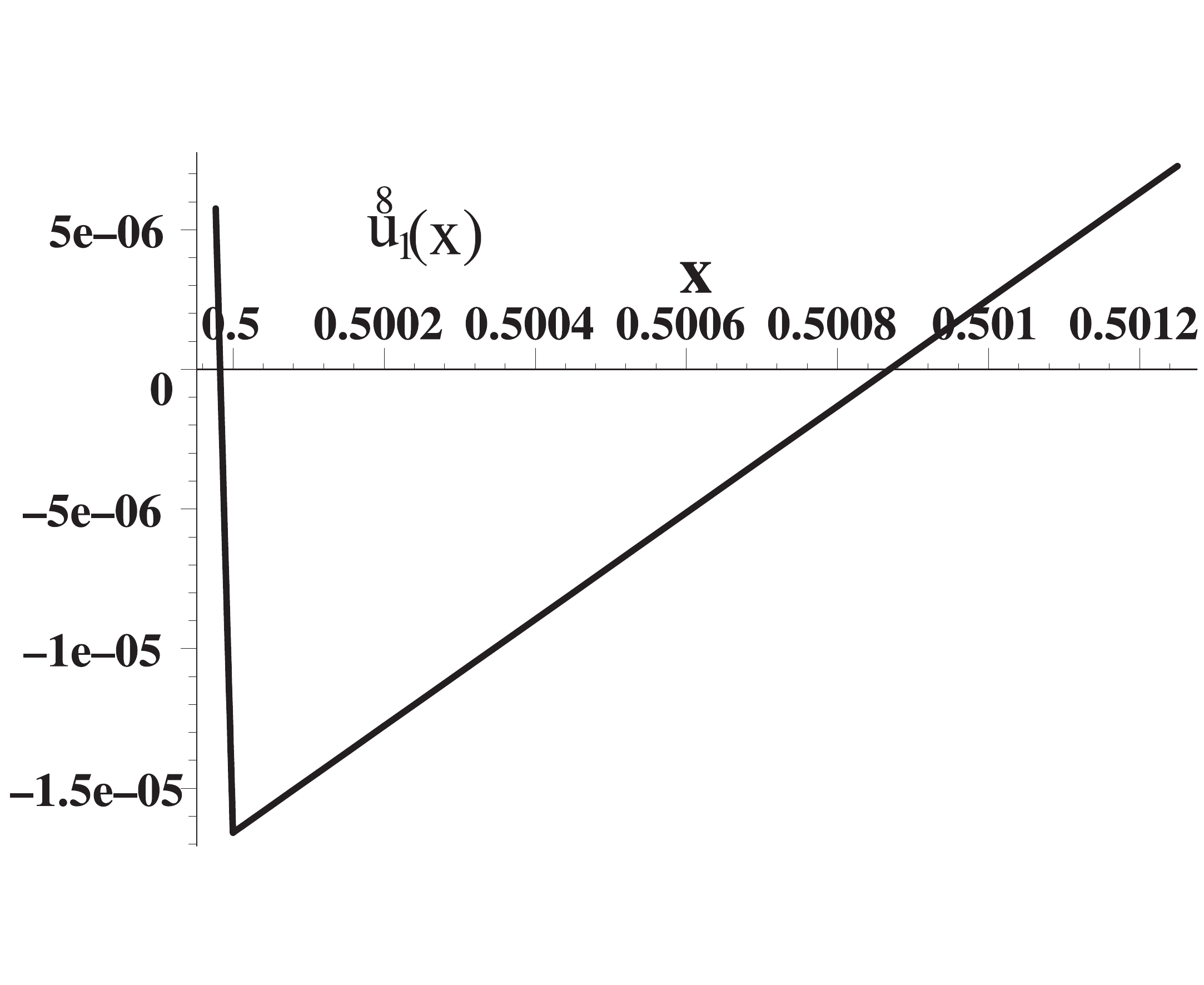}} \\ b) }
\end{minipage}
\caption{Example 4.2. FD-approximation of  $u_{1}(x)$ in $[0,1]$ (a) and $u_{1}(x)$ in $(0.5,1]$ (b).}\label{fig_7}
\end{minipage}
\end{figure}

\begin{figure}[htbp]
\begin{minipage}[h]{1\linewidth}
\center{\rotatebox{-0}{\includegraphics[width=0.48\linewidth]{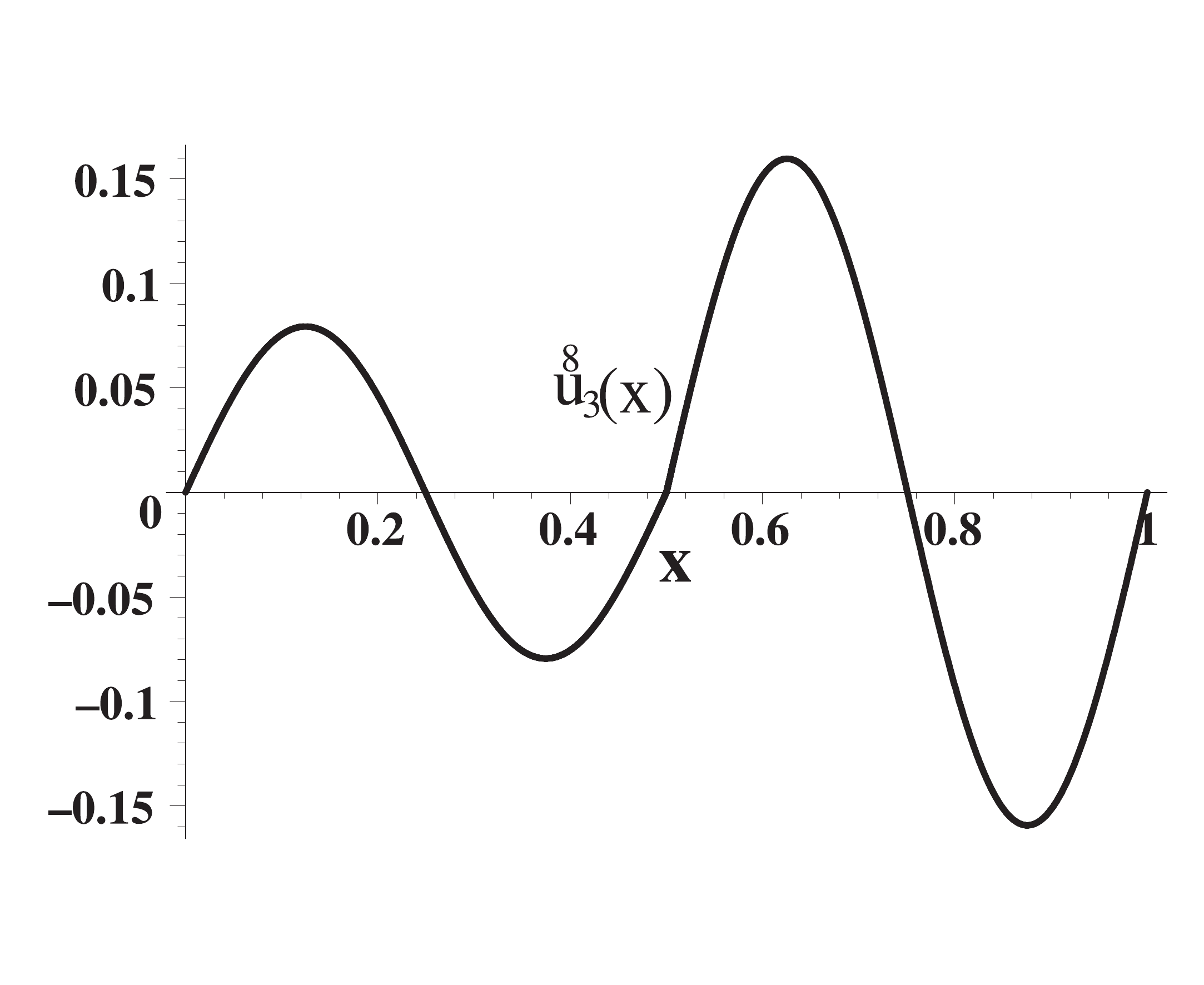}}  }
\end{minipage}
\caption{Example 4.2. FD-approximation of $u_{3}(x)$.}\label{fig_8}
\end{figure}

\begin{figure}[htbp]
\begin{minipage}[h]{1\linewidth}
\center{\rotatebox{-0}{\includegraphics[width=0.48\linewidth]{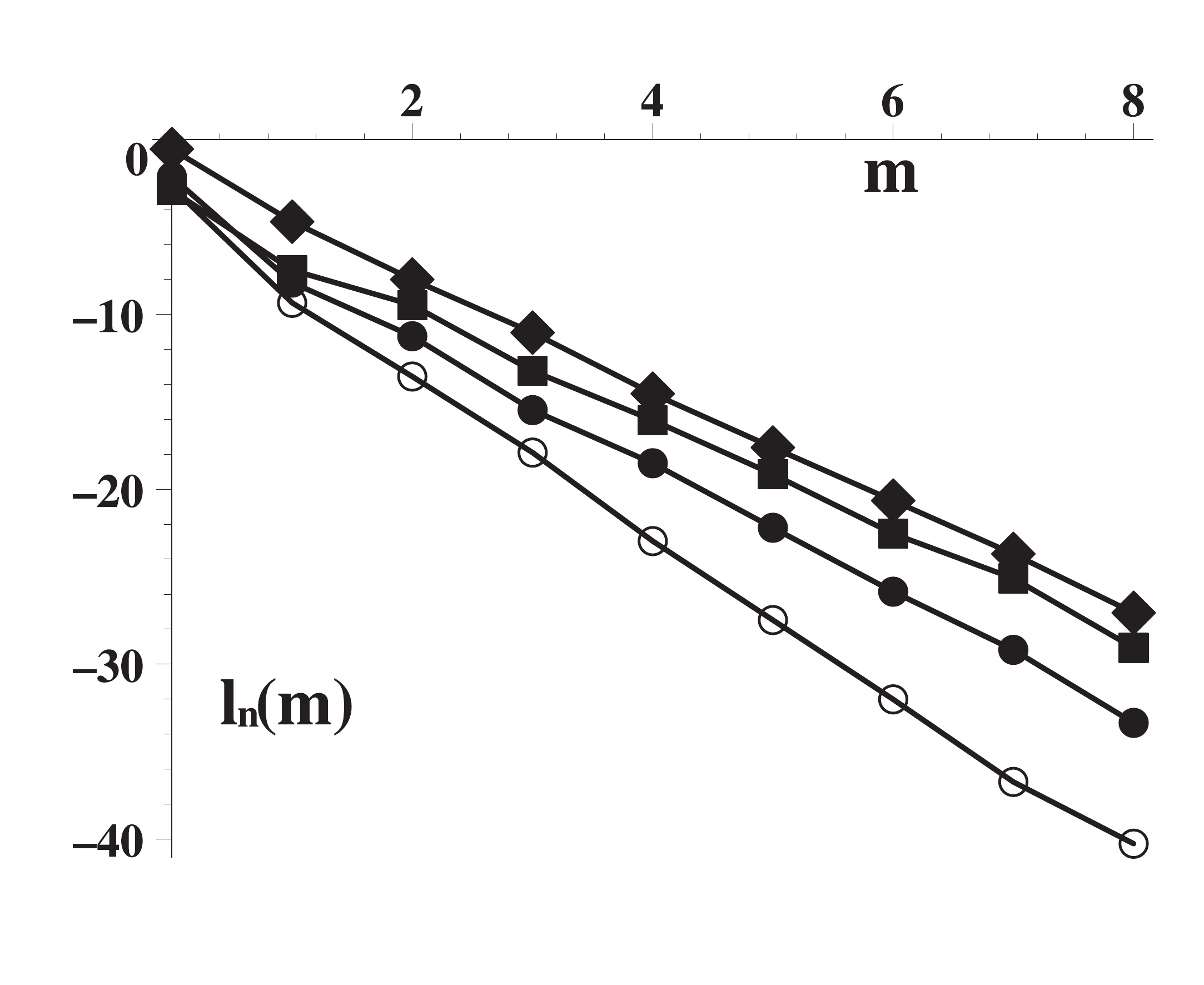}}  }
\end{minipage}
\caption{Example 4.2.
$L_\infty$-norm of $\nu_{n}^{(m)}(x)$  versus the rank $m$ of the approximation on a semi-logarithmic scale.  $\blacklozenge\colon n=0,$ $\bullet\colon n=1,$ $\blacksquare\colon n=2,$ $\circ\colon n=3.$}\label{fig_9}
\end{figure}

The results of calculations for the least four eigenvalues and corresponding eigenfunctions are presented in Tables \ref{tabl_10}-\ref{tabl_13} and Figures \ref{fig_6} - \ref{fig_9}.

\begin{table}[htbp]
\caption{Example 4.2. FD-approximation of $\lambda_{0}$ and $u_{0}(x).$}\label{tabl_10}
\begin{tabular}{|c|c|c|c|c|}
               \hline
               $m$ & $\overset{m}{\lambda}_{0}$ & $\|u_{0,1}^{(m)}(x)\|_{\infty, [0, 1/2]}$ & $\|u_{0,2}^{(m)}(x)\|_{\infty, [1/2, 1]}$ & $\|\nu_{0}^{(m)}\|_{\infty, [0, 1]}$ \\
               \hline
               0 & 17.545963379714401 & 0.24 & 0.24 & 0.58 \\
               \hline
               1 & 21.814604661812502 & 0.24e-1 & 0.24e-1 & 0.91e-2 \\
               \hline
               2 & 21.733545015000821 & 0.15e-2 & 0.15e-2 & 0.33e-3 \\
               \hline
               3 & 21.734895246851330 & 0.67e-4 & 0.67e-4 & 0.16e-4 \\
               \hline
               4 & 21.734885594786305 & 0.24e-5 & 0.24e-5 & 0.49e-6 \\
               \hline
               5 & 21.734887489687971 & 0.90e-7 & 0.90e-7 & 0.22e-7 \\
               \hline
               6 & 21.734887360755933 & 0.42e-8 & 0.42e-8 & 0.11e-8 \\
               \hline
               7 & 21.734887362837029 & 0.18e-9 & 0.18e-9 & 0.51e-10 \\
               \hline
               8 & 21.734887362829545 & 0.73e-11 & 0.73e-11 & 0.18e-11 \\
               \hline
             \end{tabular}
\newline
$\overset{8}{\lambda_{0}}= 21.734887362829545.$
\end{table}

\begin{table}[htbp]
\caption{Example 4.2. FD-approximation of $\lambda_{1}$ and $u_{1}(x).$}\label{tabl_11}
\begin{tabular}{|c|c|c|c|c|}
               \hline
               $m$ & $\overset{m}{\lambda}_{1}$ & $\|u_{1,1}^{(m)}(x)\|_{\infty, [0, 1/2]}$ & $\|u_{1,2}^{(m)}(x)\|_{\infty, [1/2, 1]}$ & $\|\nu_{1}^{(m)}\|_{\infty, [0, 1]}$ \\
               \hline
               0 & 39.478417604357434 & 0.16 & 0 & 0.12 \\
               \hline
               1 & 41.751445051880136 & 0.15e-2 & 0.31e-2 & 0.28e-3 \\
               \hline
               2 & 41.751095410924888 & 0.17e-4 & 0.26e-4 & 0.13e-4 \\
               \hline
               3 & 41.751103408518688 & 0.69e-7 & 0.20e-5 & 0.19e-6 \\
               \hline
               4 & 41.751101792373055 & 0.10e-7 & 0.30e-7 & 0.91e-8 \\
               \hline
               5 & 41.751101774628877 & 0.19e-9 & 0.12e-8 & 0.23e-9 \\
               \hline
               6 & 41.751101775581723 & 0.63e-11 & 0.35e-10 & 0.59e-11 \\
               \hline
               7 & 41.751101775606647 & 0.21e-12 & 0.62e-12 & 0.21e-12 \\
               \hline
               8 & 41.751101775606187 & 0.33e-14 & 0.33e-13 & 0.33e-14 \\
               \hline
             \end{tabular}
\newline
$\overset{8}{\lambda_{1}}= 41.751101775606187.$
\end{table}

\begin{table}[htbp]
\caption{Example 4.2. FD-approximation of $\lambda_{2}$ and $u_{2}(x).$}\label{tabl_12}
\begin{tabular}{|c|c|c|c|c|}
               \hline
               $m$ & $\overset{m}{\lambda}_{2}$ & $\|u_{2,1}^{(m)}(x)\|_{\infty, [0, 1/2]}$ & $\|u_{2,2}^{(m)}(x)\|_{\infty, [1/2, 1]}$ & $\|\nu_{2}^{(m)}\|_{\infty, [0, 1]}$ \\
               \hline
               0 & 70.183853518857661 & 0.12 & 0.12 & 0.55e-1 \\
               \hline
               1 & 73.434236059947037 & 0.36e-2 & 0.36e-2 & 0.57e-3 \\
               \hline
               2 & 73.462388394109106 & 0.39e-4 & 0.39e-4 & 0.78e-4 \\
               \hline
               3 & 73.462228765029615 & 0.33e-5 & 0.33e-5 & 0.18e-5 \\
               \hline
               4 & 73.462193970765803 & 0.92e-7 & 0.92e-7 & 0.11e-6 \\
               \hline
               5 & 73.462193827166322 & 0.52e-8 & 0.52e-8 & 0.48e-8 \\
               \hline
               6 & 73.462193886381801 & 0.21e-9 & 0.21e-9 & 0.16e-9 \\
               \hline
               7 & 73.462193887216886 & 0.89e-11 & 0.89e-11 & 0.13e-10 \\
               \hline
               8 & 73.462193887097591 & 0.54e-12 & 0.54e-12 & 0.24e-12 \\
               \hline
             \end{tabular}
\newline
$\overset{8}{\lambda_{2}}= 73.462193887097591.$
\end{table}

\begin{table}[htbp]
\caption{Example 4.2. FD-approximation of $\lambda_{3}$ and $u_{3}(x).$}\label{tabl_13}
\begin{tabular}{|c|c|c|c|c|}
               \hline
               $m$ & $\overset{m}{\lambda}_{3}$ & $\|u_{3,1}^{(m)}(x)\|_{\infty, [0, 1/2]}$ & $\|u_{3,2}^{(m)}(x)\|_{\infty, [1/2, 1]}$ & $\|\nu_{3}^{(m)}\|_{\infty, [0, 1]}$ \\
               \hline
               0 & 157.9136704174297379014 & 0.89e-1 & 0.16 & 0.66e-1 \\
               \hline
               1 & 160.2464778440802534274 & 0.50e-3 & 0.14e-2 & 0.88e-4 \\
               \hline
               2 & 160.2466449084238753618 & 0.18e-5 & 0.26e-5 & 0.13e-5 \\
               \hline
               3 & 160.2466411888247570166 & 0.33e-8 & 0.62e-7 & 0.17e-7 \\
               \hline
               4 & 160.2466412105927434146 & 0.47e-10 & 0.63e-9 & 0.11e-9 \\
               \hline
               5 & 160.2466412109016291827 & 0.50e-12 & 0.33e-11 & 0.12e-11 \\
               \hline
               6 & 160.2466412109020454979 & 0.33e-14 & 0.37e-13 & 0.12e-13 \\
               \hline
               7 & 160.2466412109020583699 & 0.31e-16 & 0.42e-15 & 0.11e-15 \\
               \hline
               8 & 160.2466412109020585522 & 0.36e-18 & 0.34e-17 & 0.33e-17 \\
               \hline
             \end{tabular}
\newline
$\overset{8}{\lambda_{3}}= 160.2466412109020585522.$
\end{table}

Similarly to Example 4.1, analyzing the results of calculations we have observed that there also exist eigenfunctions with some indexes $n$ such that their numerical approximation $\overset{8}{u}_n(x)$ have exactly $n+1$ zeros in the interval $(0,1)$ (see figure \ref{fig_7} for $n=1$). As in Example 4.1, zero approximations of such eigenpairs are determined by formula (\ref{12}) with odd $n$. Thus, taking into account the numerical calculations and convergence result, we can conclude that there does not exist an eigenfunction of problem \eqref{ex_2} with exactly $4k+1$ $(k=0,1,2,\ldots)$ zeros in the interval $(0,1)$. Moreover, there exist two linearly independent eigenfunctions with exactly $4k+2$ zeros.

Thus, Example 4.2 also demonstrates un-applicability of analytical theory developed in \cite{heinz4, zhidkov} for the case of nonself-adjoint differential operator.

The numerical experiment also confirms the exponential convergence rate of the proposed method improving along with the increase of the index of a trial eigenvalue (see figure \ref{fig_9}).

\section{Conclusions}\label{s_5}

In the paper we present a recursive algorithm for solving  both linear and nonlinear eigenvalue transmission problems with potential belonging to space $L_1(0,1)$, discontinuous flux and continuous solution. Using the method of generating functions, we find the sufficient conditions providing the exponential convergence rate of the method. Also, we establish the fact that the convergence rate of the method improves along with the increase of the index of a trial eigenpair. Numerical results confirm the theoretical conclusions. Furthermore, numerical examples and convergence results demonstrate the crucial importance of the assumption that differential operator is self-adjoint for the applicability of the theory for nonlinear eigenvalue problems developed in \cite{heinz4} and \cite{zhidkov}. As it follows from the numerical examples, the fact that differential operator is nonself-adjoint results in the lack of strong dependence between the index of an eigenpair and the numbers of zeros of the corresponding eigenfunction.


\bibliographystyle{cmam}
\bibliography{references_stat}

\end{document}